\documentclass[pdflatex,sn-mathphys-num]{sn-jnl}% Math and Physical Sciences Numbered Reference Style

\usepackage{amsmath,amsthm,amssymb,amsfonts}
\usepackage{enumerate}
\usepackage{mathtools, mathrsfs}
\usepackage[utf8]{inputenc}
\usepackage[T1]{fontenc}
\usepackage{newtxtext,newtxmath}

\usepackage{tikz}
\usetikzlibrary{calc}

\newcommand{\R}{\mathbb{R}}
\newcommand{\T}{\mathbb{T}}

\newcommand{\dist}{\textnormal{dist}}
\newcommand{\dir}{\textnormal{dir}}

\newcommand{\pb}{\textnormal{planebrush}}

\theoremstyle{thmstyleone}% numbered results: bold head, italic body
\newtheorem{theorem}{Theorem}
\newtheorem{lemma}[theorem]{Lemma}
\newtheorem{proposition}[theorem]{Proposition}
\newtheorem{corollary}[theorem]{Corollary}

\theoremstyle{thmstylethree}% definitions: bold head, upright body
\newtheorem{definition}[theorem]{Definition}
\theoremstyle{thmstyletwo}% remarks: italic head, upright body

\begin{document}

\title[A Kakeya estimate for sticky sets using a planebrush]{A Kakeya estimate for sticky sets using a planebrush}

\author*[1]{\fnm{Neeraja} \sur{Kulkarni}}\email{nkulkar7@ur.rochester.edu}

\affil[1]{\orgdiv{Department of Mathematics}, \orgname{University of Rochester}, \orgaddress{\city{Rochester}, \state{New York}, \country{USA}}}

\abstract{\unboldmath A Besicovitch set is a compact subset of $\R^n$ that contains a unit line segment in every direction. The Kakeya conjecture asserts that every Besicovitch set in $\R^n$ has Hausdorff and Minkowski dimension $n$. A recent manuscript released by OpenAI proves Sogge's local smoothing conjecture in all dimensions. Because local smoothing implies the Kakeya conjecture, that manuscript resolves the Kakeya conjecture in full. This paper gives a different proof, using incidence geometry, to a partial result for Kakeya sets in $\mathbb{R}^n$, with $n\geq 5$. It is shown that every Besicovitch set with a structural property known as \emph{stickiness} has Hausdorff dimension at least $0.60376707287\,n + \mathcal{O}(1)$. The main tool is the \emph{$k$-planebrush}, a higher-dimensional analogue of Wolff's hairbrush, which yields a Kakeya estimate for families of tubes that are $(k-1)$-plany. It is also shown, following Katz, {\L}aba and Tao, that an x-ray transform estimate at a dimension $d$ below this exponent would imply that every Besicovitch set has upper Minkowski dimension greater than $d$.}


\maketitle

\section{Introduction}
\label{sec:introduction}

\subsection{Background}
A \emph{Besicovitch set} (or \emph{Kakeya set}) is a compact set $K \subseteq \R^n$ that contains a unit line segment in every direction. When $n \geq 2$, such sets can have Lebesgue measure zero. The \emph{Kakeya set conjecture} says that every Besicovitch set in $\R^n$ has Hausdorff and Minkowski dimension $n$. The conjecture is usually studied through a discretized formulation. For $0 < \delta < 1$, a \emph{$\delta$-tube} is the $\delta$-neighborhood of a unit line segment, and the \emph{Kakeya maximal conjecture}, formulated by Bourgain~\cite{Bourgain1991}, asks for a sharp $L^p$ bound on the overlap function $\sum_{T} \chi_T$ of a family of $\delta$-tubes pointing in $\delta$-separated directions (see Definition~\ref{def-maximal-estimate} below). The maximal conjecture implies the set conjecture.

Here is a brief history of the work done on the Kakeya conjecture. The Kakeya maximal conjecture was proved in $\mathbb{R}^2$ by C\'ordoba~\cite{Cordoba1977} in 1977. In higher dimensions, Wolff~\cite{Wolff1995} introduced the ``hairbrush'' argument and proved that Besicovitch sets in $\R^n$ have Hausdorff dimension at least $(n+2)/2$. Katz, {\L}aba and Tao~\cite{Katz-Laba-Tao} studied a structural property of sets of $\delta$-tubes which they called \emph{stickiness}. They showed that a Besicovitch set of nearly minimal Minkowski dimension must be sticky, and they used this to improve the Minkowski dimension bound for Besicovitch sets in $\R^3$. Later developments include the arithmetic approach of Katz and Tao~\cite{Katz-Tao}, polynomial partitioning bounds in high dimensions due to Hickman, Rogers and Zhang~\cite{Hickman2019ImprovedBF}, and Zahl's $k$-linear Kakeya estimates~\cite{ZAHL2021107596}. Katz and Zahl~\cite{Katz-Zahl-planebrush} introduced the \emph{planebrush}, a higher-dimensional version of the hairbrush, to prove a Kakeya maximal estimate in $\R^4$. Wang and Zahl proved the Kakeya conjecture first for sticky Besicovitch sets in $\R^3$~\cite{Wang2022StickyKS} in 2022 and then for all Besicovitch sets in $\R^3$~\cite{WangZahl2025Kakeya3D} in 2025.
A recent manuscript released by OpenAI~\cite{OpenAI2026LocalSmoothing} proves Sogge's local smoothing conjecture for the wave equation in every dimension. The local smoothing conjecture is known to imply the Kakeya conjecture, so that manuscript settles the Kakeya conjecture in all dimensions.

Though the Kakeya conjecture is resolved by the recent OpenAI manuscript, this paper gives an alternative proof of Hausdorff dimension estimates for sticky Kakeya sets in $\mathbb{R}^n$ for $n\geq 5$. Unlike the proof given by the OpenAI paper, the proof here is based only on incidence geometry techniques,  including a higher-dimensional analogue of Wolff's hairbrush argument and multi-scale analysis of extremal sticky sets. 

\subsection{Notation}
Throughout, the ambient space is $\R^n$ with $n \geq 5$, and $\delta \in (0,1)$ is a small scale. The notation $A \lesssim B$ means that $A \leq CB$ for a constant $C$ depending only on $n$ (and $k$ where relevant), and $A \sim B$ means that $A \lesssim B$ and $B \lesssim A$. The notation $A \lessapprox B$ means that for every $\varepsilon > 0$ there is $C_\varepsilon$ with $A \leq C_\varepsilon \delta^{-\varepsilon} B$. In particular, factors of $\log(1/\delta)$ are absorbed by $\lessapprox$. The $r$-neighborhood of a set $A$ is denoted $\mathcal{N}_r(A)$.

\emph{Tubes.} For $0<\delta<1$, a $\delta$-tube $T$ is the $\delta$-neighborhood of a unit line segment, so $|T| \sim \delta^{n-1}$ (Subsection~\ref{subsec:sticky-kakeya} uses the normalization of Wang and Zahl, which differs only by constant factors). Its direction is $\dir(T) \in \mathbb{S}^{n-1}$. For tubes $T$, $T'$ and a linear subspace $V$, $\angle(T,T')$ denotes the angle between $\dir(T)$ and $\dir(T')$, and $\angle(T,V)$ the angle between $\dir(T)$ and $V$. Two $\delta$-tubes are \emph{essentially parallel} if $|\dir(T)-\dir(T')|\leq\delta$, and a family of $\delta$-tubes \emph{points in $\delta$-separated directions} if no two of its tubes are essentially parallel. When convenient, tubes are normalized so that $\sum_{T \in \T}|T| \sim 1$, that is, $\#\T \sim \delta^{1-n}$.

\emph{Cubes and shadings.} A \emph{$\delta$-cube} is a cube of side length $\delta$ from the standard grid. If $Z$ is a union of $r$-cubes, $\#_r Z$ denotes the number of $r$-cubes in $Z$. A \emph{shading} of a family $\T$ of $\delta$-tubes assigns to each $T \in \T$ a set $Y(T) \subseteq T$ that is a union of $\delta$-cubes; the pair is written $(\T,Y)$, or $(\T,Y)_\delta$ to indicate the scale. Write
\[E_\T := \bigcup_{T\in\T} Y(T), \qquad \T(p) := \{T \in \T : p \in Y(T)\}, \qquad \mu(p) := \#\T(p);\]
the shading is always clear from context (for example, $E_{\T'}$ refers to the shading $Y'$ of $\T'$). When no shading is specified, $Y(T) = T$. A \emph{refinement} of $(\T,Y)$ is a pair $(\T',Y')$ with $\T'\subseteq\T$, $Y'(T)\subseteq Y(T)$ for $T\in\T'$, and $\sum_{T\in\T'}|Y'(T)|\geq(\log(1/\delta))^{-C}\sum_{T\in\T}|Y(T)|$, where $C$ depends only on $n$.

\begin{definition}[Kakeya maximal estimate]
\label{def-maximal-estimate}
Let $m \geq 2$ and $1 < d \leq m$, and write $d' = d/(d-1)$. The \emph{Kakeya maximal estimate} is said to \emph{hold at dimension $d$ in $\R^m$} if for every $\varepsilon > 0$ there is $C_\varepsilon > 0$ such that
\[\Big\|\sum_{T \in \T} \chi_T\Big\|_{L^{d'}(\R^m)} \leq C_\varepsilon\, \delta^{1-\frac{m}{d}-\varepsilon}\, \big(\delta^{m-1}\#\T\big)^{1/d'}\]
for every $0<\delta<1$ and every family $\T$ of $\delta$-tubes in $\R^m$ pointing in $\delta$-separated directions.
\end{definition}

\subsection{Main results}
The first result of this paper is a Kakeya estimate for families of tubes satisfying a tangency condition called planiness.

\begin{definition}[$m$-planiness]
\label{def-planiness}
Let $C_1 \geq 1$ and $m\geq 2$. A family $(\T,Y)$ of $\delta$-tubes is called \emph{$m$-plany with constant $C_1$} if for every $\delta$-cube $Q$ there exists an $m$-plane $V(Q)$ such that
\[\angle(T,V(Q)) \leq C_1\delta \quad\text{for every } T \in \T \text{ with } Q \subseteq Y(T).\]
It is called \emph{$m$-plany} if this holds for some constant $C_1$ depending only on $n$; the constants in the estimates below then depend on $C_1$. For a point $p$, let $V(p) := V(Q)$, where $Q$ is the $\delta$-cube containing $p$.
\end{definition}
Proposition~\ref{prop-planebrush-estimate} and Theorem~\ref{prop-final-planebrush-est} below apply to $(k-1)$-plany sets, where $3\leq k \leq n-1$, and are proved using a geometric configuration called a \emph{$k$-planebrush}. This is a higher-dimensional analogue of the hairbrush of Wolff~\cite{Wolff1995}, and Katz and Zahl~\cite{Katz-Zahl-planebrush} used it to prove a Kakeya maximal estimate in $\R^4$. Roughly speaking, a $k$-planebrush at a point $p$ consists of the tubes that meet a tube through $p$ and make a fixed angle with the $(k-1)$-plane $V(p)$ supplied by the planiness hypothesis. These tubes are distributed among essentially disjoint thickened $k$-planes containing $V(p)$, so a Kakeya maximal estimate in $\R^k$ can be applied in each of them separately.
\begin{proposition}
\label{prop-planebrush-estimate}
Let $3 \leq k \leq n-1$ and $\varepsilon > 0$. Suppose that $\T$ is a family of $\delta$-tubes satisfying the following hypotheses:
\begin{enumerate}[i.]
    \item the tubes in $\T$ point in $\delta$-separated directions;
    \item each tube $T \in \T$ has a shading $Y(T)$ with $|Y(T)|\geq \lambda |T|$;
    \item the family $\T$ with shading $Y$ is $(k-1)$-plany.
\end{enumerate}
Then there exists a constant $C_\varepsilon > 0$ such that
\begin{equation}
    \label{est-planebrush}
    |E_\T| \geq C_\varepsilon\, \delta^{e(k)+\varepsilon}\, \lambda^{f(k)} \Big(\sum_{T\in \T} |T|\Big),
\end{equation}
where $d(k)$ is any exponent at which the Kakeya maximal estimate holds in $\R^k$, and
\[e(k) := \frac{(n-k)(n-2+k-d(k)) + n-2}{3(n-k)+2}, \qquad f(k) := \frac{(n-k)(d(k)+3)+3}{3(n-k)+2}.\]
\end{proposition}
Hickman, Rogers and Zhang~\cite{Hickman2019ImprovedBF} and, independently, Zahl~\cite[Theorem~1.5]{ZAHL2021107596} proved the Kakeya maximal estimate in $\R^k$ at dimension
\[\max_{2\leq \ell\leq k}\min\Big(k-\ell+2,\, \frac{k^2+k+\ell^2-\ell}{2k}\Big) \geq (2-\sqrt{2})k \qquad \text{for every } k \geq 2,\]
and Wolff~\cite{Wolff1995} proved it at dimension $(k+2)/2$. Taking
\[d(k) = d_W(k) := \max\Big\{(2-\sqrt{2})k,\, \frac{k+2}{2}\Big\}\]
in Proposition~\ref{prop-planebrush-estimate} gives the following explicit bound. The maximum equals $(k+2)/2$ for $k \leq 11$ and $(2-\sqrt{2})k$ for $k \geq 12$.

\begin{theorem}
\label{prop-final-planebrush-est}
Let $3\leq k \leq n-1$ and $\varepsilon > 0$. If $\T$ is a $(k-1)$-plany family of $\delta$-tubes pointing in $\delta$-separated directions, and each $T\in \T$ has a shading $Y(T) \subseteq T$ with $|Y(T)| \geq \lambda |T|$, then there exists a constant $C_\varepsilon > 0$ such that
\[ |E_\T| \geq C_\varepsilon\, \lambda^{f(k)}\, \delta^{n-d_{\mathrm{pl}}(k)+\varepsilon} \sum_{T\in \T} |T|,\]
where $f(k)$ is as in Proposition~\ref{prop-planebrush-estimate} with $d(k) = d_W(k)$, and
\[d_{\mathrm{pl}}(k) := n - \frac{(n-k)(n-2+k-d_W(k)) + n-2}{3(n-k)+2}.\]
For $k \geq 12$,
\[d_{\mathrm{pl}}(k) = \frac{2n+2-(\sqrt{2}-1)k}{3} - \frac{n-2-2(\sqrt{2}-1)k}{3(3(n-k)+2)}.\]
\end{theorem}

The remaining results concern sticky sets.

\begin{definition}[Sticky set]
A family $\T$ of $\delta$-tubes pointing in $\delta$-separated directions is \emph{sticky} if, for every $\varepsilon > 0$ and every $\theta$ with $\delta < \theta < 1$, there is a family $\T_\theta$ of $\theta$-tubes pointing in $\theta$-separated directions such that
\begin{enumerate}[i.]
    \item $\#\T \gtrsim \delta^{1-n+\varepsilon}$;
    \item $\#\T_\theta \gtrsim \theta^{1-n+\varepsilon}$; and
    \item $\T$ can be partitioned into the sets \[\T[T_\theta] := \{T\in \T: T\subseteq T_\theta\}, \qquad T_\theta \in \T_\theta,\]
    with
    \[\#\T[T_\theta]  \gtrsim (\delta/\theta)^{1-n+\varepsilon}\]
    for every $T_\theta \in \T_\theta$.
\end{enumerate}
The family $\T_\theta$ is called the family of \emph{$\theta$-parents} of $\T$.
\end{definition}

The main result of this paper is the following volume estimate for sticky families of tubes.

\begin{theorem}
\label{cor-sticky-estimate}
Let $d = 0.60376707287\,n + \mathcal{O}(1)$. For every $\kappa > 0$ there are $\varepsilon_0 > 0$ and $\delta_0 > 0$ such that, whenever $0 < \varepsilon \leq \varepsilon_0$, $0 < \delta \leq \delta_0$, and $\T$ is a sticky family of $\delta$-tubes with a shading $Y$ satisfying $|Y(T)| \geq \delta^\varepsilon |T|$ for all $T \in \T$,
\[|E_\T| \geq \delta^{n-d+\kappa}.\]
\end{theorem}

As a consequence (Corollary~\ref{cor-sticky-besicovitch}), every sticky Besicovitch set in $\R^n$, in the sense of Wang and Zahl~\cite{Wang2022StickyKS}, has Hausdorff dimension at least $0.60376707287\,n + \mathcal{O}(1)$. Section~\ref{sec:xray} records an x-ray transform estimate at every dimension below $(2-\sqrt2)n+\mathcal{O}(1)$, obtained by interpolating the Kakeya maximal estimate of Hickman, Rogers and Zhang~\cite{Hickman2019ImprovedBF} and Zahl~\cite{ZAHL2021107596} with the x-ray estimate of {\L}aba and Tao~\cite{Laba1999xray} (Proposition~\ref{prop-x-ray}). Following Katz, {\L}aba and Tao~\cite{Katz-Laba-Tao}, it shows that an x-ray estimate at any dimension $d < 0.60376707287\,n+\mathcal{O}(1)$ would give a Minkowski dimension bound strictly larger than $d$ for every Besicovitch set (Corollary~\ref{cor-Minkowski}).

\subsection{Outline}
Section~\ref{sec:planebrush} proves the planebrush estimates, Proposition~\ref{prop-planebrush-estimate} and Theorem~\ref{prop-final-planebrush-est}, after a two-ends reduction and a broad--narrow reduction. Section~\ref{sec:sticky} develops the multi-scale structure of extremal sticky sets, proves Theorem~\ref{cor-sticky-estimate}, and relates the sticky families of tubes defined here to the sticky Besicovitch sets of Wang and Zahl. Section~\ref{sec:xray} proves the x-ray estimate and the reduction of Minkowski dimension bounds to x-ray estimates.

\section{A Kakeya estimate for plany families of tubes}
\label{sec:planebrush}
    
    \subsection{Two-ends reduction}
    The ``two-ends reduction'' was introduced by Wolff~\cite{Wolff-x-ray} and has since been used by many authors to prove Kakeya estimates.
    \begin{definition}[Two-ends condition]
    Let $\varepsilon > 0$. A set $F\subseteq \R^n$ obeys the \emph{``two-ends'' condition with exponent $\varepsilon$} if 
    \[|F\cap B(x,r)| \lesssim r^{\varepsilon} |F|\]
    for all balls $B(x,r)$.
    \end{definition}
    Here the implicit constant is allowed to be a fixed power of $\log(1/\delta)$. In particular, if $F$ obeys the two-ends condition with exponent $\varepsilon$ and $F' \subseteq F$ satisfies $|F'| \gtrsim (\log(1/\delta))^{-C}|F|$, then $F'$ also obeys it, with a larger power of $\log(1/\delta)$ as the implicit constant.

    \begin{proposition}
    \label{prop-planebrush-estimate-with-two-ends}
    Let $3 \leq k \leq n-1$ and $\varepsilon > 0$. If $\T$ is a family of $\delta$-tubes satisfying the following hypotheses:
    \begin{enumerate}[i.]
        \item The tubes in $\T$ point in $\delta$-separated directions,
        \item Each tube $T \in \T$ has a shading $Y(T)$ with $|Y(T)| \sim \lambda |T|$,
        \item The family $\T$ with shading $Y$ is ($k$-1)-plany,
        \item For every $T\in \T$, $Y(T)$ obeys the two-ends condition with exponent $\varepsilon$,
    \end{enumerate} 
    then there exists a constant $C_\varepsilon > 0$ such that
    \[|E_\T| \geq C_\varepsilon \delta^{e(k)+\varepsilon} \lambda^{f(k)} \Big(\sum_{T\in \T} |T|\Big),\]
    where $d(k)$ is any exponent at which the Kakeya maximal estimate holds in $\R^k$, and $e(k)$, $f(k)$ are as in Proposition~\ref{prop-planebrush-estimate}.
    \end{proposition}
    
    \begin{lemma} [Two-ends reduction]
    \label{lem-two-ends-reduction}
    Proposition $\ref{prop-planebrush-estimate-with-two-ends}$ implies Proposition $\ref{prop-planebrush-estimate}$.
    \end{lemma}
    The proof of Lemma~\ref{lem-two-ends-reduction} is standard; see~\cite{Tao1999Lectures}. The argument there carries over to the present setting because the only rescaling it uses is isotropic: a portion of each tube of length $\sim r$ is mapped to a $(\delta/r)$-tube by $x \mapsto x/r$. An isotropic rescaling preserves angles between tubes, and between tubes and planes. Hence, if $\angle(T, V(Q)) \lesssim \delta$ before rescaling, the same bound holds afterwards. Since $\delta \leq \delta/r$, the rescaled family of $(\delta/r)$-tubes is again $(k-1)$-plany, so hypothesis (iii) is preserved. Finally, a dyadic pigeonholing step passes to a subfamily on which $|Y(T)| \sim \lambda'|T|$ for some $\lambda' \gtrsim \lambda$, as required in hypothesis (ii) of Proposition~\ref{prop-planebrush-estimate-with-two-ends}. This costs only logarithmic factors, which are absorbed into $\delta^{-\varepsilon}$ exactly as at the end of the proof of Lemma~\ref{lem-broad-narrow-reduction}, and it preserves the two-ends condition with a larger power of $\log(1/\delta)$ as the implicit constant.
    
    \subsection{Broad-narrow reduction}
    If $\T$ is a family of $\delta$-tubes with shading $Y$, a point $p\in E_\T$ is called \emph{$s$-narrow} if there exists a unit vector $v_p\in\mathbb{S}^{n-1}$ such that 
    \[\#\{T\in \T(p): \angle(\dir(T),v_p) \leq s\} \geq \frac{1}{2} \mu(p).\]
    If this condition fails at a point $p\in E_\T$, then $p$ is called \emph{$s$-broad}. 
    
    \begin{proposition}
    \label{prop-conditional-planebrush-est}
    Let $3 \leq k \leq n-1$ and $\varepsilon > 0$. If $\T$ is a family of $\delta$-tubes satisfying the following hypotheses:
    \begin{enumerate}[i.]
        \item The tubes in $\T$ point in $\delta$-separated directions,
        \item Each tube $T \in \T$ has a shading $Y(T)$ with $|Y(T)| \sim \lambda |T|$,
        \item The family $\T$ with shading $Y$ is ($k$-1)-plany,
        \item For every $T\in \T$, $Y(T)$ obeys the two-ends condition with exponent $\varepsilon$,
        \item There is a set $G \subseteq E_\T$ satisfying the condition
        \[\sum_{T\in \T} |Y(T) \cap G| \gtrsim \frac{1}{\log(1/\delta)} \sum_{T\in \T} |Y(T)|, \]
        for which any point $p \in G$ satisfies:
        \[\#\{(T_1, T_2) \in \T(p) \times \T(p): \angle(T_1, T_2) \gtrsim \delta^\varepsilon\} \gtrsim \frac{1}{\log(1/\delta)}\mu(p)^2; \]
    \end{enumerate} 
    then there exists a constant $C_\varepsilon > 0$ such that
    \[|E_\T| \geq C_\varepsilon \delta^{e(k)+\varepsilon} \lambda^{f(k)} \Big(\sum_{T\in \T} |T|\Big)\]
    where $d(k)$ is any exponent at which the Kakeya maximal estimate holds in $\R^k$, and $e(k)$, $f(k)$ are as in Proposition~\ref{prop-planebrush-estimate}.
    \end{proposition}
    \begin{lemma}[Broad-narrow reduction]
    \label{lem-broad-narrow-reduction}
    Proposition $\ref{prop-conditional-planebrush-est}$ implies Proposition $\ref{prop-planebrush-estimate-with-two-ends}$.        
    \end{lemma}
    The proof of the broad-narrow reduction relies on the following simple observations about broad and narrow points.
    \begin{lemma} \label{lem-s-broad-angle-separation}
        If a point $p \in E_\T$ is $s$-broad, then 
        \begin{equation}\label{angle-sep-condition}
        \#\bigg\{(T_1, T_2) \in \T(p)\times \T(p): \angle(T_1, T_2) > s\bigg\} \geq \frac{3}{4}\, \mu(p)^2. 
        \end{equation}
    \end{lemma}
    \begin{proof}
        If a point $p$ is $s$-broad with respect to $(\T,Y)$, then one may find a set of $4s$-separated directions $\{v_i\in \mathbb{S}^{n-1}\}_{i=1}^{Cs^{1-n}}$, where $C\sim 1$ is a fixed constant, such that for every $i$,
    \[\#\{T\in \T(p): \angle(\dir(T),v_i) \leq s\} \leq \frac{1}{2} \mu(p).\]
    Since the set of directions $\{v_i\in \mathbb{S}^{n-1}\}_{i=1}^{Cs^{1-n}}$ is $4s$-separated, it follows that if $T_i \in \{T\in \T: \angle(\dir(T), v_i)\leq s\}$, $T_j \in \{T\in \T: \angle(\dir(T), v_j)\leq s\}$, and $i\neq j$, then
    \[\angle(T_i, T_j) \geq s.\]
    The total number of ordered pairs $(T_1, T_2) \in \T(p) \times \T(p)$ in which $T_1$ and $T_2$ both lie in the $s$-neighborhood of the same direction is less than $\frac{1}{4} \mu(p)^2$, and the conclusion follows.
    \end{proof}
    \begin{lemma}
    \label{lem-exists-broad-and-narrow-set}
    If $\T$ is a family of $\delta$-tubes with shading $Y$, there is a number $s$ and a set $G \subseteq E_\T$ such that 
     \[\sum_{T\in \T} |Y(T) \cap G| \gtrsim \frac{1}{\log(1/\delta)} \sum_{T\in \T} |Y(T)| \]
    and every point $p \in G$ satisfies two conditions:
    \begin{enumerate}[(i)]
        \item $p$ is $s$-narrow, and
        \item the following estimate holds:
        \[\#\bigg\{(T_1,T_2)\in \T(p)\times \T(p): \angle(T_1,T_2) \geq \frac{s}{4}\bigg\} \geq \frac{3}{4} \mu(p)^2.\]
    \end{enumerate}
    \end{lemma}
    \begin{proof}
        For every point $p$, define $k_p$ by
        \begin{equation} \label{def_k_p}
        k_p = \max\bigg\{k\in \mathbb{Z} : \exists v_p \in \mathbb{S}^{n-1} \text{ such that } \#\{T\in \T(p): \angle(v_p, \dir(T)) \lesssim 2^{-k} \} \geq \frac{1}{2} \mu(p)\bigg\}.\end{equation}
        This maximum exists because $k$ must lie between 0 and $\sim \log(\delta^{-1})$. By dyadic pigeonholing, one may choose a value of $s$ such that
        \begin{equation}\label{def_pigeonholing_s}
        \frac{s}{2}\leq 2^{-k_p} < s\end{equation}
        for every point $p$ in a set $G \subseteq E_\T$, such that
        \[\sum_{T\in \T} |Y(T) \cap G| \gtrsim \frac{1}{\log(1/\delta)} \sum_{T\in \T} |Y(T)|. \]
        Every point $p \in G$ satisfies the two required conditions.
        To see this, note from $(\ref{def_k_p})$ and $(\ref{def_pigeonholing_s})$ that $p$ satisfies the first condition. For the second condition, the definition of $k_p$ implies that at every $p \in G$, 
        \[\forall v \in \mathbb{S}^{n-1}: \#\bigg\{T\in \T(p): \angle(\dir(T),v) \geq \frac{s}{4} \bigg\} \geq \frac{1}{2} \mu(p),\]
        which means $p$ is $s/4$-broad.
        Applying Lemma $\ref{lem-s-broad-angle-separation}$ gives
        \[\# \bigg\{(T_1, T_2) \in \T(p) \times \T(p): \angle(T_1, T_2) > \frac{s}{4} \bigg\} \geq \frac{3}{4} \mu(p)^2.\]
    \end{proof}
    
    \begin{proof}[Proof of Lemma \ref{lem-broad-narrow-reduction}]
        Assume that Proposition $\ref{prop-conditional-planebrush-est}$ is true. Let $\varepsilon > 0$ be given and let $\T$ be a family of $\delta$-tubes that satisfies hypotheses (i)--(iv) of Proposition $\ref{prop-planebrush-estimate-with-two-ends}$.
        Apply Lemma $\ref{lem-exists-broad-and-narrow-set}$ and let $G$ and $s$ be as in the statement of that lemma.
        For every point $p \in G$, let $v_p\in \mathbb{S}^{n-1}$ be as in $(\ref{def_k_p})$.
        Define a new shading by
        \[Y'(T) = \{p \in Y(T)\cap G: \angle(v_p, \dir(T))\leq s\},\]
        and correspondingly define
        \[\T'(p) = \{T\in \T: p\in Y'(T)\}.\]
        By Lemma $\ref{lem-exists-broad-and-narrow-set}$, for every $p \in G$, 
        \[\#\{T\in \T: p\in Y'(T)\} \geq \frac{1}{2} \mu(p),\]
        and hence
        \[\sum_{T\in \T} |Y'(T)| \sim \int_G \#\T'(p)\; dp \gtrsim \int_G \mu(p)\; dp \sim \sum_{T\in \T}|Y(T)\cap G|.\]
        Also from Lemma $\ref{lem-exists-broad-and-narrow-set}$,
        \[\sum_{T\in \T} |Y(T)\cap G| \gtrsim \frac{1}{\log(1/\delta)} \sum_{T\in \T}|Y(T)|.\]
        Combining the above two statements gives
        \begin{equation}
        \label{eq-shading-comparable-1}
        \sum_{T\in \T} |Y'(T)| \gtrsim \frac{1}{\log(1/\delta)}\sum_{T\in \T}|Y(T)|.
        \end{equation}
        Next, by dyadic pigeonholing, choose a subfamily $\T'' \subseteq \T$, a shading $Y''$ of $\T''$ with $Y''(T) \subseteq Y'(T)$ for every $T \in \T''$, and a number $0 < \lambda' < 1$ such that
        \[|Y''(T)| \sim \lambda' |T| \text{ for all } T\in \T''\]
        and in addition
        \begin{equation}
        \label{eq-shading-comparable-2}
        \sum_{T\in \T''} |Y''(T)| \gtrsim \frac{1}{\log(1/\delta)}\sum_{T\in \T}|Y'(T)|.
        \end{equation}
        From $(\ref{eq-shading-comparable-1})$, $(\ref{eq-shading-comparable-2})$, hypothesis (ii) of Proposition $\ref{prop-planebrush-estimate-with-two-ends}$ and $\#\T'' \leq \#\T$, it follows that
        \begin{equation}
        \label{eq-shading-comparable-3}
        \lambda' \gtrsim \frac{\lambda}{\left(\log(1/\delta)\right)^2}.
        \end{equation}
        It is possible to cover the upper hemisphere of $\mathbb{S}^{n-1}$ with $\sim s^{1-n}$ finitely overlapping tiles of side length $s$, where a tile of side length  $s$ means the image via the map $\xi \mapsto (\xi, \sqrt{1-|\xi|^2})$ of translates of tiles $[0,s)^{n-1}$ tiling $\R^{n-1}$. For each tile $\tau$ of side length $s$, define 
        \[\T(\tau) := \{T\in \T'': \dir(T) \in \tau\}.\]
        Let $\tau_1$ and $\tau_2$ be two non-adjacent tiles so $\dist(\tau_1, \tau_2) > 2s$. 
        If $T_1 \in \T(\tau_1)$ and $T_2 \in \T(\tau_2)$, the angle between the tubes would be greater than $2s$:
        \[\angle(T_1, T_2) > 2s.\]
        By the definition of the shading $Y''$, this means that $Y''(T_1) \cap Y''(T_2) = \emptyset$. Thus, under the shading $Y''$, the tubes in $\T(\tau_1)$ and those in $\T(\tau_2)$ have disjoint shadings.
        Within each $s$-tile $\tau$ covering $\mathbb{S}^{n-1}$, there are at most $\lesssim (s/\delta)^{n-1}$ tubes of $\T$.  Without loss of generality, assume that the center of the tile $\tau$ is the north pole of the sphere, i.e. the point $(0,\ldots, 0, 1)$. Perform a non-isotropic rescaling $\phi$ which takes $(x_1,\ldots, x_n) \mapsto (x_1/s, \ldots, x_{n-1}/s, x_n)$. After this rescaling,  $\phi(T)$ is contained in the $n$-fold dilate of a $\delta/s$-tube $\Tilde{T}$ that points in the same direction as $T$. Define a shading $\Tilde{Y}$ for each $\Tilde{T}$ by \[\Tilde{Y}(\Tilde{T}):= \phi(Y''(T)).\]
        By $(\ref{eq-shading-comparable-3})$,
        \[|\Tilde{Y}(\Tilde{T})| \sim \lambda' |\Tilde{T}|.\]
        It remains to check that $(\Tilde{\T}, \Tilde{Y})$ satisfies hypotheses (i) --- (v) in the statement of Proposition $\ref{prop-conditional-planebrush-est}$.
        Hypothesis (ii) has already been shown, and hypothesis (i) holds because of the non-isotropic way in which the re-scaling $\phi$ was done. The family $\Tilde{\T}$ is $(k-1)$-plany because the family $\T$ is assumed to be $(k-1)$-plany. For hypothesis (iv), observe that $Y''(T) \subseteq Y(T)$ and, by hypothesis (ii) and $(\ref{eq-shading-comparable-3})$, $|Y''(T)| \sim \lambda'|T| \gtrsim (\log(1/\delta))^{-2}|Y(T)|$. Since $Y(T)$ obeys the two-ends condition with exponent $\varepsilon$, so does $Y''(T)$, with a larger power of $\log(1/\delta)$ as the implicit constant, and hence so does $\Tilde{Y}(\Tilde{T})$.
        To show that $(\Tilde{\T}, \Tilde{Y})$ satisfies hypothesis (v), note that at every point $p \in \bigcup_{T\in \T(\tau)} Y'(T)$, Lemma $\ref{lem-exists-broad-and-narrow-set}$ gives
        \[\#\bigg\{(T_1,T_2)\in \T'(p)\times \T'(p): \angle(T_1, T_2)\gtrsim \frac{s}{4}\bigg\} \geq \frac{3}{4}(\#\T'(p))^2.\]
        Finally, after the rescaling, 
        \[\angle(T_1, T_2)\geq s/4 \text{ implies that }\angle(\Tilde{T}_1, \Tilde{T}_2) \gtrsim 1,\]
        since $\phi$ enlarges angles between directions in the tile $\tau$ by a factor of $\sim s^{-1}$. 
        Thus the family $(\Tilde{\T}, \Tilde{Y})$ satisfies hypotheses (i)--(v) of Proposition $\ref{prop-conditional-planebrush-est}$ at scale $\delta/s$, with density $\lambda'$ in hypothesis (ii) and with $\varepsilon/2$ in place of $\varepsilon$. (The two-ends condition with exponent $\varepsilon$ implies the two-ends condition with exponent $\varepsilon/2$, and the angle bound in (v) holds with a constant independent of $\delta$.) Write
        \[a := e(k) \quad\text{and}\quad b := f(k),\]
        where $e(k)$ and $f(k)$ are as in Proposition~\ref{prop-planebrush-estimate}. Note that $a > 0$ and $b \geq 1$.
        By Proposition $\ref{prop-conditional-planebrush-est}$, there exists $C_\varepsilon > 0$ such that
        \[|E_{\Tilde{\T}}| \geq C_\varepsilon \left(\frac{\delta}{s}\right)^{a+\varepsilon/2} (\lambda')^{b} \sum_{\Tilde{T} \in \Tilde{\T}} |\Tilde{T}|.\]
        Undoing the rescaling $\phi$ multiplies both $|E_{\Tilde{\T}}|$ and $\sum |\Tilde{T}|$ by $s^{n-1}$. Since $s < 1$ and $a + \varepsilon/2 > 0$, the factor $s^{-a-\varepsilon/2}$ is at least $1$, and therefore
        \[\left|\bigcup_{T\in \T(\tau)} Y''(T) \right|\geq C_\varepsilon\, \delta^{a+\varepsilon/2} (\lambda')^{b} \sum_{T\in \T(\tau)} |T|.\]
        The tiles can be divided into $C_n$ classes, each consisting of pairwise non-adjacent tiles, where $C_n$ depends only on $n$. Choose a class $\mathcal{C}$ for which
        \[\sum_{\tau \in \mathcal{C}} \sum_{T\in \T(\tau)} |T| \geq C_n^{-1} \sum_{T\in \T''} |T|.\]
        For distinct tiles in $\mathcal{C}$ the sets $\bigcup_{T\in\T(\tau)} Y''(T)$ are disjoint, so summing the previous estimate over $\tau \in \mathcal{C}$ gives
        \begin{equation}
        \label{eq-broad-narrow-lambda-prime}
        |E_\T| \geq |E_{\T''}| \geq C_n^{-1} C_\varepsilon\, \delta^{a+\varepsilon/2} (\lambda')^{b} \sum_{T\in \T''} |T|.
        \end{equation}
        It remains to replace $\lambda'$ by $\lambda$ and $\T''$ by $\T$, at the cost of powers of $\log(1/\delta)$. By $(\ref{eq-shading-comparable-1})$, $(\ref{eq-shading-comparable-2})$ and hypothesis (ii),
        \[\lambda' \sum_{T\in\T''} |T| \sim \sum_{T\in\T''} |Y''(T)| \gtrsim \frac{1}{(\log(1/\delta))^2} \sum_{T\in\T} |Y(T)| \geq \frac{\lambda}{(\log(1/\delta))^2} \sum_{T\in\T} |T|,\]
        and by $(\ref{eq-shading-comparable-3})$,
        \[(\lambda')^{b-1} \gtrsim \frac{\lambda^{b-1}}{(\log(1/\delta))^{2(b-1)}}.\]
        Multiplying these two estimates gives
        \[(\lambda')^{b} \sum_{T\in \T''} |T| \gtrsim \frac{\lambda^{b}}{(\log(1/\delta))^{2b}} \sum_{T\in \T} |T|.\]
        Since $(\log(1/\delta))^{2b} \leq C'_\varepsilon\, \delta^{-\varepsilon/2}$ for some constant $C'_\varepsilon$, combining this with $(\ref{eq-broad-narrow-lambda-prime})$ yields
        \[|E_\T| \geq C''_\varepsilon\, \delta^{e(k)+\varepsilon} \lambda^{f(k)} \sum_{T\in \T} |T|,\]
        which is the conclusion of Proposition $\ref{prop-planebrush-estimate-with-two-ends}$. This completes the proof.
    \end{proof}
    
    \subsection{Volume of a planebrush}
    \label{subsec:planebrush-volume}
    The proof of Proposition~\ref{prop-planebrush-estimate} rests on the $k$-planebrush, which is modelled on Wolff's hairbrush. Wolff's hairbrush in $\R^3$ consists of the $\delta$-tubes that meet a fixed central $\delta$-tube $T$ at large angles. Each such tube lies in the $\delta$-neighborhood of a plane containing the axis of $T$. These thickened planes are essentially disjoint away from $T$, and each tube of the hairbrush lies in exactly one of them. Applying the two-dimensional Kakeya maximal estimate separately in each thickened plane and summing gives a lower bound for the volume of the hairbrush.

    \begin{definition}[$k$-planebrush]
    \label{def-planebrush} Let $(\T,Y)$ be $(k-1)$-plany. If $p \in E_\T$ and $0 < \rho < 1$, the $k$-\emph{planebrush at the point $p$ with angle $\rho$} is the set of tubes
    \[\T_{\pb}(p,\rho) = \{T\in \T: \angle(T,V(p)) \sim \rho \text{ and }\exists T' \in \T(p): T\cap T' \neq \emptyset \}.\]
    \end{definition}

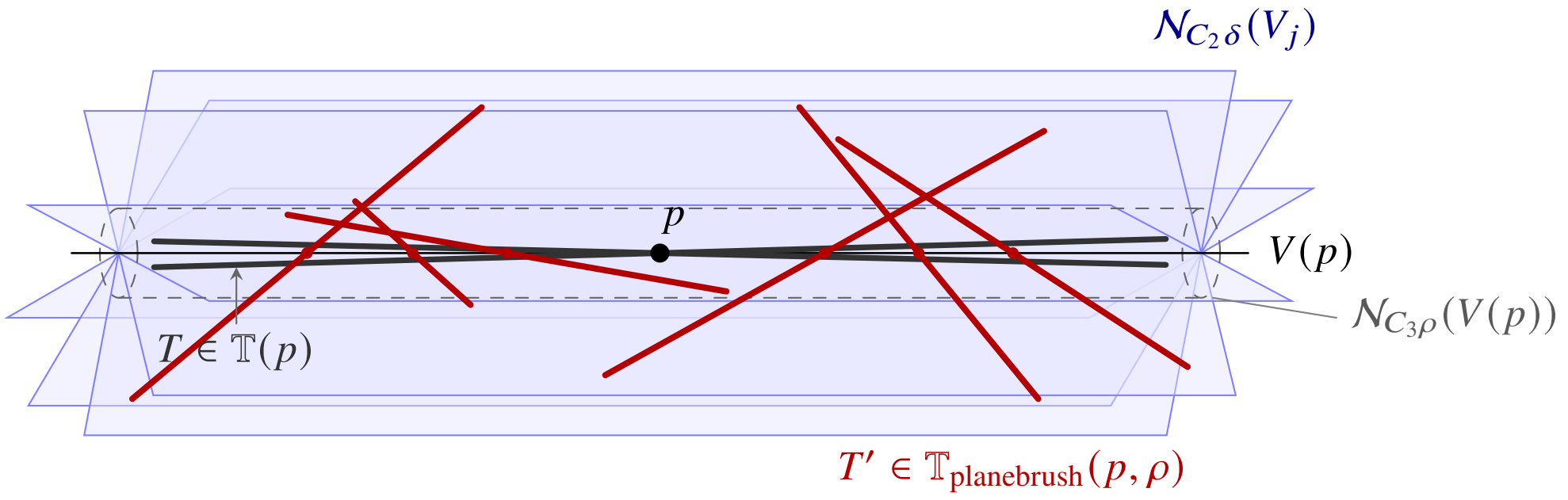
\begin{figure}[t]
\centering
\resizebox{\linewidth}{!}{%
\begin{tikzpicture}[scale=1.0,line cap=round,
  page/.style={fill=blue!10,draw=blue!50,line width=0.4pt,fill opacity=0.45},
  brush/.style={line width=1.5pt,red!70!black},
  core/.style={line width=1.4pt,black!80}]
  \def\L{4.6}
  % u(a) = cos(a)*(0.95,0.55) + sin(a)*(0,1.45): projected unit vector in V(p)^perp
  % thickened k-planes V_j = span(V(p), alpha_j), all containing V(p)
  \foreach \a in {0,36,72,108,144} {
    \coordinate (u) at ({0.95*cos(\a)},{0.55*cos(\a)+1.45*sin(\a)});
    \filldraw[page] ($(-\L,0)-(u)$) -- ($(\L,0)-(u)$) -- ($(\L,0)+(u)$) -- ($(-\L,0)+(u)$) -- cycle;
  }
  % rho-neighborhood of V(p) (removed in the volume estimate)
  \draw[dashed,black!60] (-\L,0) ellipse (0.16 and 0.38);
  \draw[dashed,black!60] (-\L,0.38) -- (\L,0.38);
  \draw[dashed,black!60] (-\L,-0.38) -- (\L,-0.38);
  \draw[dashed,black!60] (\L,0) ellipse (0.16 and 0.38);
  % the (k-1)-plane V(p)
  \draw[black,line width=0.5pt] (-\L-0.4,0) -- (\L+0.4,0);
  % tubes through p, tangent to V(p)
  \draw[core] (-4.3,-0.12) -- (4.3,0.12);
  \draw[core] (-4.3,0.10) -- (4.3,-0.10);
  % planebrush tubes: each meets a tube through p and lies in exactly one page
  \foreach \x/\a/\s in {-3.0/72/1, -1.3/144/-1, 1.4/36/1, 3.0/108/-1, -2.1/0/-1, 2.2/72/-1} {
    \coordinate (u) at ({0.95*cos(\a)},{0.55*cos(\a)+1.45*sin(\a)});
    \coordinate (q) at (\x,0);
    \draw[brush] ($(q)-(1.25,0)-\s*0.8*(u)$) -- ($(q)+(1.25,0)+\s*0.8*(u)$);
    \fill[red!70!black] (q) circle (1.4pt);
  }
  \fill (0,0) circle (2.3pt);
  \node[anchor=south,inner sep=1pt] at (0.12,0.12) {$p$};
  % labels
  \node[anchor=west] at (\L+0.45,0) {$V(p)$};
  \node[anchor=west,black!65] at (\L+1.15,-0.55) {$\mathcal{N}_{C_3\rho}(V(p))$};
  \draw[black!50] (\L+1.15,-0.55) -- (\L+0.1,-0.38);
  \node[anchor=south,blue!55!black] at ({\L+0.95*cos(72)},{0.55*cos(72)+1.45*sin(72)+0.05}) {$\mathcal{N}_{C_2\delta}(V_j)$};
  \node[anchor=north,red!70!black] at (3.0,-1.55) {$T'\in\T_{\pb}(p,\rho)$};
  \node[anchor=north,black!80] at (-3.6,-0.55) {$T\in\T(p)$};
  \draw[black!60,->,>=stealth] (-3.6,-0.6) -- (-3.6,-0.12);
\end{tikzpicture}%
}
\caption{A schematic $k$-planebrush at a point $p$. The horizontal line represents the $(k-1)$-plane $V(p)$. The black tubes pass through $p$ and are tangent to $V(p)$. The shaded sheets are the thickened $k$-planes $\mathcal{N}_{C_2\delta}(V_j)$, where $V_j=\mathrm{span}(V(p),\alpha_j)$. They all contain $V(p)$ and are essentially disjoint outside the dashed neighborhood $\mathcal{N}_{C_3\rho}(V(p))$. Each red tube belongs to the planebrush $\T_{\pb}(p,\rho)$: it makes angle $\sim\rho$ with $V(p)$, meets a tube through $p$, and lies in exactly one sheet.}
\label{fig:planebrush}
\end{figure}

    In a $k$-planebrush, the central tube is replaced by a $(k-1)$-plane. At a point $p$, the planiness hypothesis provides a $(k-1)$-plane $V(p)$ to which every tube through $p$ is tangent. The planebrush $\T_{\pb}(p,\rho)$ of Definition~\ref{def-planebrush} consists of the tubes that make angle $\sim\rho$ with $V(p)$ and meet one of the tubes through $p$. Each of these tubes lies in exactly one of $\lesssim (\rho/\delta)^{n-k}$ $\delta$-thickened $k$-planes containing $V(p)$. Away from a $\rho$-neighborhood of $V(p)$, these thickened $k$-planes are essentially disjoint, so the best available Kakeya maximal estimate in $\R^k$ can be applied separately inside each of them (Fig.~\ref{fig:planebrush}). The $(k-1)$-planiness hypothesis guarantees that a planebrush contains enough tubes for this estimate to be effective.

    This subsection estimates the volume of a $k$-planebrush. The first step is to estimate the number of tubes it contains.
    
    \begin{lemma}[Constant multiplicity refinement]
    \label{lem-const-mult-refinement}
        Let $\T$ be a family of $\delta$-tubes pointing in $\delta$-separated directions. Let $Y$ be a shading of $\T$, and let $0< \lambda < 1$ be such that $|Y(T)| \geq \lambda |T|$ for every $T\in \T$. Then there exists a refinement $\T'\subset \T$ and a new shading  $Y'(T) \subseteq Y(T)$ for each $T\in \T'$, which satisfy the following properties:
        \begin{enumerate}[i]
            \item There is a dyadic number $\mu$ such that for every point $p \in E_{\T'}$, $ \mu \leq \#\T'(p) < 2\mu$.
            \item $\sum_{T\in\T'} |Y'(T)| \gtrsim \left(\sum_{T\in\T} |Y(T)| \right)/\log(1/\delta)^2$.
            \item $\mu \gtrapprox \lambda \delta^{n-1} \#\T/|E_\T|$.
            \item There is a number $\lambda' \gtrapprox \lambda$ such that $|Y'(T)| \sim \lambda'|T|$ for every $T \in \T'$.
        \end{enumerate}
    \end{lemma}
    \begin{proof}
        For every integer $j$ between $1$ and $\log \delta^{1-n}$, let $A_j \subseteq  E_\T$ be the set of points with multiplicity around $2^j$, i.e.
        \[A_j = \bigg\{ p \in E_\T : 2^j \leq \#\T(p) < 2^{j+1} \bigg\}.\]
        For each $j$, define a shading $Y_j$ by $Y_j(T) := Y(T) \cap A_j.$
        For each integer $\ell$ between 0 and $\log (\lambda/\delta)$, let $\T_{j,\ell}$ be the collection of tubes $T$ for which the shading $Y_j(T)$ has density around $2^{-\ell} \lambda $, i.e.
        \[\T_{j,\ell} = \bigg\{T\in \T: 2^{-\ell} \lambda |T| \leq |Y_j(T)| <  2^{-\ell+1} \lambda |T| \bigg\}. \]
        From the definitions of $j$ and $\ell$, one can observe that
        \[\sum_{T\in \T} |Y(T)| \sim \#\T\cdot \lambda\delta^{n-1} \sim \sum_{j=0}^{\log(1/\delta)} \sum_{\ell = 0}^{\log(1/\delta)} \#\T_{j,\ell} 2^{-\ell} \lambda \delta^{n-1}.  \]
        By the pigeonhole principle, there exists at least one pair $(j,\ell)$ that satisfies
        \begin{equation}
        \label{eq-high-incidence-pair}
        \#\T_{j,\ell} 2^{-\ell} \lambda \delta^{n-1} \gtrsim \frac{1}{(\log(1/\delta))^2} \#\T \cdot \lambda \delta^{n-1}.
        \end{equation}
        Take $\mu = 2^j$, define the shading $Y'$ by $Y_j$ and define $\T'$ as $\T_{j,\ell}$ for this choice of $j$ and $\ell$.
        From here, points (i) and (ii) in the statement of the lemma follow immediately from $(\ref{eq-high-incidence-pair})$.
        One can also observe that
        \[\mu = 2^j \geq \frac{\lambda \cdot \#\T_{j,\ell}\cdot 2^{-\ell}\delta^{n-1}}{|E_\T|} \gtrapprox \frac{\lambda\,\delta^{n-1}\#\T}{|E_\T|}, \]
        which is statement (iii).
        Finally, every $T \in \T' = \T_{j,\ell}$ satisfies $|Y'(T)| \sim 2^{-\ell}\lambda|T|$, and since $\#\T_{j,\ell} \leq \#\T$, $(\ref{eq-high-incidence-pair})$ gives $2^{-\ell} \gtrsim (\log(1/\delta))^{-2}$. This is statement (iv) with $\lambda' := 2^{-\ell}\lambda$.
    \end{proof}

    \begin{lemma}[Planebrush multiplicity]
    \label{lem-planebrush-mult}
    Let $\T$ be a family of direction-separated $\delta$-tubes with shading $Y$ that obeys the hypotheses of Proposition $\ref{prop-conditional-planebrush-est}$, and write $S := \sum_{T\in\T}|T|$. Then there exist an angle $\delta \leq \rho \leq 1$ and a point $p$ such that
    \[\#\T_\pb(p,\rho) \gtrapprox \frac{\lambda^3 S^2 \,\rho \,\delta^{-1}}{|E_\T|^2}.\]
    \end{lemma}
    \begin{proof}
    Let $\T', Y'$ be the constant multiplicity refinement of $\T,Y$ given by Lemma $\ref{lem-const-mult-refinement}$. By Lemma~\ref{lem-const-mult-refinement}(i), every point $p \in E_{\T'}$ satisfies $\#\T'(p) \geq \mu$; fix such a point $p$.
    For this point $p$, one can count the number of triples $(T_1, Q_1, T_2)$ where $T_1$ is a tube whose shading contains $p$ and $T_2$ is a tube whose shading intersects the shading of $T_1$ in a $\delta$-cube $Q_1$. More precisely,
    \begin{align*}
    \#\{(T_1, Q_1, T_2)\in \T'\times \mathcal{Q}\times \T: T_1 \in \T'(p), Q_1 \subseteq Y'(T_1) \cap Y'(T_2) \} &\gtrsim \mu \cdot \lambda' \delta^{-1} \cdot \mu  \\ & \gtrapprox \mu^2 \lambda \delta^{-1}, 
    \end{align*}
    where $\lambda' \gtrapprox \lambda$ is the density from Lemma~\ref{lem-const-mult-refinement}(iv). By the pigeonhole principle, there exists an angle $\rho$ with $\delta \leq \rho \leq 1$ such that
     \begin{align*}
     \#\{(T_1, Q_1, T_2)\in \T'\times \mathcal{Q}\times \T: T_1 \in \T'(p), Q_1 \subseteq Y'(T_1) \cap Y'(T_2), \angle(T_2, V(p)) \sim \rho\} &\\ 
     \gtrapprox \mu^2 \lambda \delta^{-1}.& \end{align*}
     It remains to show that each pair $(Q_1, T_2)$ arises from only $\lessapprox 1$ tubes $T_1$. By hypothesis (ii) of Proposition~\ref{prop-conditional-planebrush-est}, $|Y(T_1)| \sim \lambda|T_1|$, and by Lemma~\ref{lem-const-mult-refinement}(iv), $|Y'(T_1)| \sim \lambda'|T_1| \gtrsim (\log(1/\delta))^{-2}|Y(T_1)|$. So by hypothesis (iv), for every $T_1 \in \T'(p)$ and every $r > 0$,
     \[|Y'(T_1) \cap B(p,r)| \leq |Y(T_1)\cap B(p,r)| \lesssim r^{\varepsilon}(\log(1/\delta))^{C}\,|Y'(T_1)|\]
     for a constant $C$. Choose $r_0 := c\,(\log(1/\delta))^{-C/\varepsilon}$, with $c > 0$ small enough that the right-hand side is at most $\frac{1}{2}|Y'(T_1)|$ when $r = r_0$. Then at least half of each $Y'(T_1)$ lies outside $B(p, r_0)$, so restricting the triples above to cubes $Q_1$ at distance $\geq r_0$ from $p$ loses at most a constant factor. For such a cube, every $T_1 \in \T'(p)$ containing $Q_1$ points in a direction within $\lesssim \delta/r_0$ of the direction from $p$ to $Q_1$. Since the tubes point in $\delta$-separated directions, there are $\lesssim r_0^{1-n} = (\log(1/\delta))^{O(1/\varepsilon)} \lessapprox 1$ such tubes $T_1$. Therefore
     \[\#\{(Q_1, T_2)\in \mathcal{Q}\times \T: T_1 \in \T'(p), (T_1, Q_1, T_2) \text { is a triple in the previous step }\} \gtrapprox \mu^2 \lambda \delta^{-1}. \]
    Since the tube $T_2$ makes angle $\rho$ with $V(p)$, there are $\sim \rho^{-1}$ $\delta$-cubes contained in $T_2$ that are also contained in a $\delta$-neighborhood of $V(p)$. 
    Thus, for each $\delta$-tube $T_2$ that makes angle $\rho$ with $V(p)$, there are potentially up to $\rho^{-1}$ cubes $Q_1$ for which $(Q_1, T_2)$ is a pair counted in the previous statement, i.e. a pair for which $Q_1 \subseteq Y'(T_2) \cap \mathcal{N}_\delta(V(p))$ and $Q_1 \subseteq Y'(T_1)$ for some $T_1\in\T(p)$ . It follows that
    \begin{align*}\#\{T_2\in \T:  \exists T_1\in \T(p) \text{ such that } Y'(T_1) \cap Y'(T_2) \neq \emptyset \text{ and } \angle(T_2, V(p)) &\sim \rho\}\\
    & \gtrapprox \mu^2 \lambda \delta^{-1} \rho. \end{align*}
    The left-hand side is exactly $\#\T_{\pb}(p,\rho)$. Since $\mu \gtrapprox \lambda S/|E_\T|$ by Lemma~\ref{lem-const-mult-refinement}(iii), the lemma follows.
    \end{proof}

    \begin{lemma}[Volume of a single planebrush]
        \label{lem-vol-planebrush}
        Let $(\T, Y)$ satisfy the hypotheses of Proposition~\ref{prop-conditional-planebrush-est}, and let $\T_\pb(p,\rho)$ be a planebrush at a point $p$ with angle $\rho$. Then
        \[|E_{\T_\pb(p,\rho)}| \gtrapprox \delta^{n-1+k-d(k)} \lambda^{d(k)}\, \#\T_\pb(p,\rho). \]
    \end{lemma} 
    \begin{proof}
    Choose a maximal set of $\delta/\rho$-separated points $\{\alpha_j\}_j$ in the unit sphere in the space $V(p)^\perp = \{v\in \R^n: \langle v, V(p)\rangle = 0\}$, which is an $(n-k+1)$-dimensional space. For each $\alpha_j$, define a $k$-plane $V_j$ by $V_j = \textnormal{span}(V(p), \alpha_j)$. The resulting collection of $k$-planes $\{V_j\}_j$ satisfies the following two conditions.
    \begin{enumerate}[(A)]
    \item 
    Every tube $T \in \T_\pb(p,\rho)$ is contained in $\mathcal{N}_{C_2\delta}V_j$, with $C_2 \sim 1$, for some $j$.
    \item 
    If a point $x \in \R^n$ satisfies $\dist(x, V(p)) \geq C_3 \rho$, then $x \in \mathcal{N}_{C_2\delta}V_j$ for at most $\sim (C_2/C_3)^{n-k+1}$ choices of $k$-plane $V_j$.
    \end{enumerate}
    To show that (A) holds, let $T \in \T_\pb(p, \rho)$ and write $\dir(T) = v^{||} + v^\perp$, where $v^{||} \in V(p)$ and $v^{\perp} \in \left(V(p)\right)^\perp$. 
    Choose $\alpha_j$ such that $\angle(\alpha_j, v^\perp) \leq \delta/\rho$. 
    This is possible because of the maximality of the set $\{\alpha_j\}$. 
    Write $v^{\perp} = v_1 + v_2$, where $v_1$ lies in the direction of $\alpha_j$ and $\langle v_2, \alpha_j\rangle = 0$. 
    Since $\angle(T, V(p)) \lesssim \rho$, it follows that $|v^\perp| = \sin \left(\angle(T, V(p))\right) \lesssim \rho$. 
    Now 
    \[\sin\left(\angle(T, V_j)\right) = |v_2| = \sin\left(\angle(v_2, v^{\perp})\right)\cdot|v^\perp| \lesssim \delta, \]
    so (A) is proved.
    
    To show (B), let a point $x$ satisfy $\dist(x, V(p)) \geq C_3 \rho$. 
    If $x = x^{||} + x^\perp$, where $x^{||}$ lies in $V(p)$ and $\langle x^\perp, x^{||}\rangle = 0$, then $|x^\perp| \gtrsim C_3 \rho$. 
    Suppose that $x \in \mathcal{N}_{C_2 \delta} V_j$ for some $j$. Then 
    \begin{equation}
    \label{cond-B}
    \dist\left(\frac{x^\perp}{|x^\perp|}, \alpha_j\right) \lesssim \frac{C_2\delta}{C_3\rho}. 
    \end{equation}
    Since the directions $\{\alpha_j\}$ are $\delta/\rho$-separated, there are only at most $\sim (C_2/C_3)^{n-k+1}$ choices of $\alpha_j$ that can satisfy $(\ref{cond-B})$. This proves (B).
    
    For each $k$-plane $V_j$, let $\T_j$ be the family of tubes given by \[\T_j = \{T\in \T_\pb(p,\rho): T \subset \mathcal{N}_{C_2\delta}V_j \}.\]
    Define the set $F_j$ by
    \[F_j = E_{\T_j} \setminus \mathcal{N}_{C_3 \rho} V(p).\]
    Assume first that $\rho \geq K\delta$, where $K := C/C_3$ for a suitable constant $C$; the case $\rho < K\delta$ is treated at the end of the proof. Translate so that $p = 0$. Each $T \in \T_j$ meets a tube through $p$, which lies in $\mathcal{N}_{C\delta}(V(p))$, and makes angle $\sim \rho$ with $V(p)$. Hence $T \cap \mathcal{N}_{C_3\rho}(V(p))$ is contained in a segment of $T$ of length $\lesssim C_3 + \delta/\rho \lesssim C_3$, and so in a ball $B(x_T, C'C_3)$. By hypothesis (iv) of Proposition~\ref{prop-conditional-planebrush-est} (the two-ends condition, with its fixed exponent $\varepsilon$),
    \[|Y(T) \cap \mathcal{N}_{C_3\rho}(V(p))| \leq |Y(T) \cap B(x_T, C'C_3)| \lesssim (C'C_3)^{\varepsilon}(\log(1/\delta))^{C}\,|Y(T)|.\]
    Choose $C_3 := c\,(\log(1/\delta))^{-C/\varepsilon}$ with $c > 0$ small enough that the right-hand side is at most $\frac{1}{2}|Y(T)|$. Then
    \[|Y(T) \cap F_j| \geq \tfrac{1}{2}|Y(T)| \gtrsim \lambda |T|.\]
    With this choice, the overlap bound in (B) is $(C_2/C_3)^{n-k+1} = (\log(1/\delta))^{O(1/\varepsilon)} \lessapprox 1$, and $K \lessapprox 1$.
    Summing over all $T\in\T_j$ gives
    \begin{equation} 
    \label{eq-L_1}
    \int_{F_j} \sum_{T\in\T_j} \chi_{Y(T)} \geq \frac{\lambda}{2} \delta^{n-1} \#\T_j.
    \end{equation}
    On the other hand, H\"older's inequality gives
    \begin{equation} 
    \label{eq-holder}
    \int_{F_j} \sum_{T\in\T_j} \chi_{Y(T)} \leq || \sum_{T\in\T_j} \chi_T ||_{L^{d(k)'}(\R^n)} \cdot |F_j|^{1/d(k)},
    \end{equation}
    where $d(k)$ is the exponent in Proposition~\ref{prop-conditional-planebrush-est}, at which the Kakeya maximal estimate holds in $\R^k$,
    and $d(k)'$ is the H\"older conjugate of $d(k)$.
    A Kakeya maximal estimate will now be applied to $|| \sum_{T\in\T_j} \chi_T ||_{L^{d(k)'}}$. 
    Let $\pi_{V_j}$ be the orthogonal projection from $\R^n$ onto $V_j$. 
    Since each $T \in \T_j$ satisfies $T \subset \mathcal{N}_{C_2\delta} V_j$, it follows that $\angle(T, V_j) \lesssim C_2\delta$. 
    This implies that $\pi_{V_j}(T)$ is a cylinder of radius $\delta$ and length $\gtrsim C_2^{-1}$ and is contained in a $\sim C_2^{-1}$-dilation of a $\delta$-tube in $\R^k$. 
    Thus,
    \begin{align*} 
    || \sum_{T\in\T_j} \chi_T ||_{L^{d(k)'}(\R^n)}
     &=  \left(\int_{\R^{n-k}}  \int_{\R^k} |\sum_{T\in\T_j} \chi_T |^{d(k)'}\right)^{1/d(k)'} \\
     &=  \left(\int_{[-C_2\delta, C_2\delta]^{n-k}}  \int_{V_j} |\sum_{T\in\T_j} \chi_{T} |^{d(k)'}\right)^{1/d(k)'} \\
     & \lesssim \delta^{(n-k)/d(k)'}\left(  \int_{V_j} |\sum_{T\in\T_j} \chi_{C_2\pi_{V_j}(T)} |^{d(k)'}\right)^{1/d(k)'}.
    \end{align*}
    After applying the Kakeya maximal estimate, this gives
    \[
    || \sum_{T\in\T_j} \chi_T ||_{L^{d(k)'}(\R^n)} \lessapprox \delta^{(n-k)/d(k)'} \delta^{-k/d(k)+1-\varepsilon} (\delta^{k-1}\#\T_j)^{1/d(k)'}.
    \]
    Combining this estimate with $(\ref{eq-L_1})$ and $(\ref{eq-holder})$ gives
    \[
    \frac{\lambda}{2} \delta^{n-1} \#\T_j \lessapprox \delta^{(n-k)/d(k)'} \delta^{-k/d(k)+1-\varepsilon} (\delta^{k-1}\#\T_j)^{1/d(k)'} \cdot |F_j|^{1/d(k)},
    \]
    which simplifies to
    \begin{equation}
        \label{F_j-estimate}
        |F_j| \gtrapprox \delta^{n-1} \lambda^{d(k)}\delta^{k-d(k)} \#\T_j.
    \end{equation}
    Summing $(\ref{F_j-estimate})$ over all $j$ and using properties (A) and (B) above gives
    \begin{equation}
        \label{E-estimate-1}
        |E_{\T_\pb(p,\rho)}| \gtrapprox \delta^{n-1} \lambda^{d(k)}\delta^{k-d(k)} \#\T_\pb(p,\rho).
    \end{equation}
    This is the claimed estimate when $\rho \geq K\delta$.

    Finally, suppose that $\rho < K\delta$. Every $T \in \T_\pb(p,\rho)$ meets a tube through $p$ and makes angle $\lesssim K\delta$ with $V(p)$, so $T \subseteq \mathcal{N}_{CK\delta}(V(p))$. In particular, all these tubes lie in $\mathcal{N}_{CK\delta}(V_1)$ for any single $k$-plane $V_1 \supseteq V(p)$. Running the argument above with the single plane $V_1$ in place of $\{V_j\}$, with $C_2$ replaced by $CK$, and with $F_1 := E_{\T_\pb(p,\rho)}$ (no neighborhood of $V(p)$ is removed, so (B) is not needed), gives the same estimate with an extra factor $K^{O(1)} \lessapprox 1$.
    \end{proof}

    \begin{lemma}[Size of a planebrush]
    \label{lem-planebrush-size}
    Let $3 \leq k \leq n-1$, let $\T$ be a family of $\delta$-tubes pointing in $\delta$-separated directions with a $(k-1)$-plany shading $Y$, and let $\delta \leq \rho \leq 1$. Then for every point $p$,
    \[\#\T_\pb(p,\rho) \lesssim \rho^{n-k+1}\delta^{1-n}.\]
    \end{lemma}
    \begin{proof}
    Every $T \in \T_\pb(p,\rho)$ satisfies $\angle(T, V(p)) \lesssim \rho$. So $\dir(T)$ lies in the $C\rho$-neighborhood, in $\mathbb{S}^{n-1}$, of the great sphere $\mathbb{S}^{n-1} \cap V(p)$, which has dimension $k-2$. This neighborhood has measure $\lesssim \rho^{(n-1)-(k-2)} = \rho^{n-k+1}$. The directions of the tubes in $\T_\pb(p,\rho)$ are $\delta$-separated and $\rho \geq \delta$, so there are $\lesssim \rho^{n-k+1}\delta^{1-n}$ of them.
    \end{proof}
    
    \subsection{Proof of \textit{k}-planebrush estimate}
    \label{subsec:planebrush-proof}
    The goal of this subsection is to prove Proposition~$\ref{prop-conditional-planebrush-est}$. A single planebrush suffices. Lemma~\ref{lem-planebrush-mult} provides a planebrush containing many tubes, Lemma~\ref{lem-planebrush-size} shows that such a planebrush must have a large angle, and Lemma~\ref{lem-vol-planebrush} converts the number of tubes into volume.

    \begin{proof}[Proof of Proposition $\ref{prop-conditional-planebrush-est}$]
        Let $\varepsilon > 0$ and let $(\T, Y)$ satisfy the hypotheses of the proposition. Write $m := n-k \geq 1$ and
        \[S := \sum_{T\in\T}|T| = \delta^{n-1}\#\T \lesssim 1,\]
        where $S \lesssim 1$ because the tubes point in $\delta$-separated directions. By Lemma~\ref{lem-planebrush-mult}, there are a point $p$ and an angle $\delta \leq \rho \leq 1$ such that
        \begin{equation}
        \label{eq-pb-1}
        \#\T_\pb(p,\rho) \gtrapprox \frac{\lambda^3 S^2 \rho\, \delta^{-1}}{|E_\T|^2}.
        \end{equation}
        By Lemma~\ref{lem-planebrush-size},
        \begin{equation}
        \label{eq-pb-2}
        \#\T_\pb(p,\rho) \lesssim \rho^{m+1}\delta^{1-n}.
        \end{equation}
        By Lemma~\ref{lem-vol-planebrush}, and since the union of the shadings of the planebrush is contained in $E_\T$,
        \begin{equation}
        \label{eq-pb-3}
        |E_\T| \gtrapprox \delta^{n-1+k-d(k)} \lambda^{d(k)}\, \#\T_\pb(p,\rho).
        \end{equation}
        Combining $(\ref{eq-pb-1})$ and $(\ref{eq-pb-2})$ gives a lower bound for the angle:
        \begin{equation}
        \label{eq-pb-angle}
        \rho^{m} \gtrapprox \frac{\lambda^3 S^2 \delta^{n-2}}{|E_\T|^2}.
        \end{equation}
        Combining $(\ref{eq-pb-1})$ and $(\ref{eq-pb-3})$ gives
        \[|E_\T|^3 \gtrapprox \delta^{n-2+k-d(k)} \lambda^{d(k)+3} S^2 \rho.\]
        Raising this to the power $m$ and inserting $(\ref{eq-pb-angle})$ yields
        \[|E_\T|^{3m+2} \gtrapprox \delta^{m(n-2+k-d(k)) + n-2}\, \lambda^{m(d(k)+3)+3}\, S^{2m+2},\]
        and taking $(3m+2)$-th roots,
        \[|E_\T| \gtrapprox \delta^{e(k)} \lambda^{f(k)} S^{(2m+2)/(3m+2)}.\]
        Since $S \lesssim 1$ and $(2m+2)/(3m+2) \leq 1$, we have $S^{(2m+2)/(3m+2)} \gtrsim S$. Unwinding the notation $\gtrapprox$, this gives
        \[|E_\T| \geq C_\varepsilon\, \delta^{e(k)+\varepsilon} \lambda^{f(k)} \sum_{T\in\T}|T|,\]
        which completes the proof.
    \end{proof}
    \begin{proof}[Proof of Proposition $\ref{prop-planebrush-estimate}$]
        This follows by combining Proposition $\ref{prop-conditional-planebrush-est}$ with Lemmas $\ref{lem-two-ends-reduction}$ and $\ref{lem-broad-narrow-reduction}$.
    \end{proof}
    \begin{proof}[Proof of Theorem $\ref{prop-final-planebrush-est}$]
        Let $\varepsilon > 0$ and let $\T, Y$ be as stated in the hypothesis of the proposition. The exponent $d(k) = d_W(k)$ is admissible in Proposition $\ref{prop-planebrush-estimate}$: the Kakeya maximal estimate holds in $\R^k$ at dimension $(2-\sqrt{2})k$ by~\cite{Hickman2019ImprovedBF} and~\cite[Theorem~1.5]{ZAHL2021107596}, and at dimension $(k+2)/2$ by~\cite{Wolff1995}. Proposition $\ref{prop-planebrush-estimate}$ therefore yields
        \[|E_\T| \geq C_\varepsilon \delta^{e(k) + \varepsilon} \lambda^{f(k)} \sum_{T\in \T} |T|,\]
        which is the claimed estimate, since $n - e(k) = d_{\mathrm{pl}}(k)$ by definition. For $k \geq 12$ we have $d_W(k) = (2-\sqrt{2})k$, and then
        \[n - e(k) = n - \frac{n-2+(\sqrt{2}-1)k}{3} - \frac{n-2-2(\sqrt{2}-1)k}{3(3(n-k)+2)},\]
        which gives the explicit formula for $d_{\mathrm{pl}}(k)$.
    \end{proof}

    \section{A Kakeya estimate for sticky sets}
    \label{sec:sticky}

    Subsection~\ref{subsec:sticky-kakeya} proves Proposition~\ref{prop-sticky-estimate}, and hence Theorem~\ref{cor-sticky-estimate}, using the multiscale structure of nearly extremal families of tubes established by Wang and Zahl~\cite{Wang2022StickyKS}. Subsection~\ref{subsec:wang-zahl} explains how sticky families of $\delta$-tubes, as defined here, relate to the sticky Besicovitch sets of Wang and Zahl.

    \subsection{Proof of the sticky estimate}
    \label{subsec:sticky-kakeya}

    This subsection proves Proposition~\ref{prop-sticky-estimate} using the multiscale structure of nearly extremal families of tubes established by Wang and Zahl~\cite[Sections~2--3]{Wang2022StickyKS}. Their main structural result says that such a family is coarsely self-similar: if the $\delta$-tubes are covered by $\rho$-tubes, where $\rho = \delta^{\tau}$, then the $\rho$-tubes themselves form a nearly extremal family. The argument applies a ``plany-or-broad'' dichotomy to the $\rho$-tubes. If they are mostly $(k-1)$-plany at scale $\rho$, the planiness is at the scale of the tubes, so Theorem~\ref{prop-final-planebrush-est} applies at scale $\rho$ and bounds the extremal exponent directly. If they are mostly $k$-broad, the broadness passes down to the $\delta$-tubes at angle $\sim\rho$, and Zahl's $k$-linear estimate is applied at scale $\delta$ with transversality $\gtrsim\rho^{k-1} = \delta^{\tau(k-1)}$, which costs only a factor $\delta^{O(\tau)}$.
    \\ \\
    \noindent\textbf{Lines, tubes and coverings.} Every line $\ell$ in $\R^n$ that is not parallel to the hyperplane $\{x_n = 0\}$ can be written as $\ell = (\textbf{a}_\ell,0)+\R e_\ell$, with $\textbf{a}_\ell\in\R^{n-1}$ and $e_\ell\in\mathbb{S}^{n-1}$; this line is also denoted $(\textbf{a}_\ell, e_\ell)$. The distance between two such lines is
    \[\dist(\ell,\ell') := |\textbf{a}_\ell-\textbf{a}_{\ell'}|+\angle(e_\ell,e_{\ell'}),\]
    and $m_{\mathcal{L}}$ denotes the product of Lebesgue measure on $\R^{n-1}$ and normalized surface measure $\nu$ on $\mathbb{S}^{n-1}$, regarded as a measure on the set $\mathcal{L}$ of such lines. 
    %With these definitions, the packing dimension of a set of lines agrees with its upper modified box dimension. 
    Following~\cite[Section~2.1]{Wang2022StickyKS}, in this section a $\delta$-tube is the $2n\delta$-neighborhood, intersected with $[-1,1]^n$, of a line $\ell = (\textbf{a}_\ell,e_\ell)$ (its \emph{coaxial line}) with $\textbf{a}_\ell\in[-\frac1n,\frac1n]^{n-1}$ and $\angle(e_\ell,e_n)\leq\pi/3$. Two $\delta$-tubes are \emph{essentially distinct} if their coaxial lines have distance greater than $\delta$. A $\rho$-tube $\tilde T$ \emph{covers} a $\delta$-tube $T$, where $\delta\leq\rho$, if their coaxial lines have distance at most $\rho/2$; in particular $\angle(T,\tilde T)\leq\rho/2$. If $\tilde\T$ covers $\T$, then $\T[\tilde T]$ denotes the tubes of $\T$ covered by $\tilde T$, and a family $(\tilde\T,\tilde Y)_\rho$ \emph{covers} $(\T,Y)_\delta$ if $\tilde\T$ covers $\T$ and $Y(T)\subseteq\tilde Y(\tilde T)$ whenever $T\in\T[\tilde T]$.

    \begin{definition}[{\cite[Definition~2.1]{Wang2022StickyKS}}] \label{def:epsilon-extremal}
        For $s,t,\delta\in(0,1)$, let $M(s,t,\delta) := \inf|E_\T|$, where the infimum is taken over all pairs $(\T,Y)_\delta$ such that
        \begin{enumerate}[(a)]
            \item the tubes in $\T$ are essentially distinct;
            \item for each $\delta\leq\rho\leq 1$, $\T$ can be covered by a set of $\rho$-tubes, at most $\delta^{-t}$ of which are essentially parallel to a common $\rho$-tube; and
            \item $\sum_{T\in\T}|Y(T)|\geq\delta^s$.
        \end{enumerate}
        Let $N(s,t) := \limsup_{\delta\to 0}\frac{\log M(s,t,\delta)}{\log\delta}$, and let $\sigma_n := \lim_{(s,t)\to(0,0)}N(s,t) = \inf_{(s,t)\in(0,1)^2}N(s,t)$. A pair $(\T,Y)_\delta$ is called \emph{$\varepsilon$-extremal} if it satisfies (a), (b) and (c) with $s=t=\varepsilon$, and $|E_\T|\leq\delta^{\sigma_n-\varepsilon}$.
    \end{definition}
    The limit defining $\sigma_n$ exists because $N$ is nondecreasing in each variable and $0\leq N\leq n$. Note that $\sigma_n$ is a codimension: the sticky Kakeya conjecture corresponds to $\sigma_n = 0$. By~\cite[Lemma~3.1]{Wang2022StickyKS}, for every $\varepsilon,\delta_1>0$ there is an $\varepsilon$-extremal pair $(\T,Y)_\delta$ for some $\delta\in(0,\delta_1]$.

    \begin{proposition}[Self-similarity of extremal families; {\cite[Proposition~3.2]{Wang2022StickyKS}}]
        \label{prop-wz-self-similar}
        For every $\varepsilon>0$ there are $\eta>0$ and $\delta_0>0$ such that the following holds for every $\delta\in(0,\delta_0]$. Let $(\T,Y)_\delta$ be $\eta$-extremal and let $\rho\in[\delta^{1-\varepsilon},\delta^{\varepsilon}]$. Then there are a refinement $(\T',Y')_\delta$ of $(\T,Y)_\delta$ and a family $(\tilde\T,\tilde Y)_\rho$ covering it, with the following properties.
        \begin{enumerate}[(i)]
            \item Every $T\in\T'$ is covered by exactly one $\tilde T\in\tilde\T$.
            \item $(\tilde\T,\tilde Y)_\rho$ is $\varepsilon$-extremal. In particular, $\sum_{\tilde T}|\tilde Y(\tilde T)|\geq\rho^{\varepsilon}$, at most $\rho^{-\varepsilon}$ tubes of $\tilde\T$ are essentially parallel to a common $\rho$-tube, and $|E_{\tilde\T}|\leq\rho^{\sigma_n-\varepsilon}$.
            \item $\#\tilde\T(p)\leq\rho^{-\sigma_n-\varepsilon}$ for every $p\in\R^n$.
            \item There is $\mu_{\mathrm{fine}}\leq(\rho/\delta)^{\sigma_n+\varepsilon}$ such that $\#\T'[\tilde T](p) \sim \mu_{\mathrm{fine}}$ for every $\tilde T\in\tilde\T$ and every $p\in E_{\T'[\tilde T]}$.
            \item There is $w>0$ such that $|Q\cap E_{\T'[\tilde T]}| = w$ for every $\tilde T\in\tilde\T$ and every $\rho$-cube $Q\subseteq\tilde Y(\tilde T)$.
        \end{enumerate}
    \end{proposition}
    Items (ii), (iii) and the bound on $\mu_{\mathrm{fine}}$ are items (i), (iii) and (iv) of~\cite[Proposition~3.2]{Wang2022StickyKS}. Item (i), the constancy of $\mu_{\mathrm{fine}}$, and item (v) are established in Steps~1 and~3 of its proof (the properties of $\T_2$ and of $(\T_5,Y_5)_\delta$, $(\tilde\T_5,\tilde Y_5)_\rho$ there), and they are preserved by the final step, which only removes some tubes of $\tilde\T_5$ together with the $\delta$-tubes they cover.

    The plany/broad dichotomy is applied with a threshold that may be small.

    \begin{lemma}[Plany or broad]
        \label{lem-k-linear-or-k-plany}
        Let $0<\rho<1$, $C_1\geq 1$, $3 \leq k \leq n-1$, and $0<a\leq\frac{1}{2k}$, and let $\mathcal{S}$ be a family of $\rho$-tubes with shading $Z$, where each $Z(S)$ is a union of $\rho$-cubes. For a $\rho$-cube $Q$, let $\mathcal{S}(Q) := \{S\in\mathcal{S}: Q\subseteq Z(S)\}$. Call $Q$ \emph{plany} if there is a $(k-1)$-plane $V(Q)$ such that
        \[\#\{S\in\mathcal{S}(Q): \angle(S,V(Q)) \leq C_1\rho\} \geq a\,\#\mathcal{S}(Q),\]
        and \emph{broad} otherwise. (Note that the definitions of plany and broad depend on $k$, $\rho$, $C_1$ and $a$.) Then one of the following holds.
        \begin{enumerate}[(a)]
            \item The shading $Z'(S) := \bigcup\{Q \subseteq Z(S): Q \text{ plany},\ \angle(S,V(Q))\leq C_1\rho\}$ satisfies $\sum_{S}|Z'(S)| \geq \frac{a}{2}\sum_S |Z(S)|$, and $(\mathcal{S}, Z')$ is $(k-1)$-plany with constant $C_1$, in the sense of Definition~\ref{def-planiness} at scale $\rho$.
            \item The shading $Z''(S) := \bigcup\{Q\subseteq Z(S): Q \text{ broad}\}$ satisfies $\sum_S |Z''(S)| \geq \frac{1}{2}\sum_S|Z(S)|$.
        \end{enumerate}
    \end{lemma}
    \begin{proof}
        Since $\sum_S |Z(S)| = \sum_Q \#\mathcal{S}(Q)\,|Q|$, either the plany cubes or the broad cubes contribute at least half of this sum. In the second case (b) holds. In the first case, 
        \begin{align*}\sum_S |Z'(S)| &= \sum_{Q \text{ plany}} \#\{S\in\mathcal{S}(Q): \angle(S,V(Q))\leq C_1\rho\}\,|Q| \\ &\geq a\sum_{Q\text{ plany}}\#\mathcal{S}(Q)\,|Q| \geq \frac{a}{2}\sum_S|Z(S)|,\end{align*}
        and every $S$ with $Q\subseteq Z'(S)$ satisfies $\angle(S,V(Q))\leq C_1\rho$ by construction, so (a) holds.
    \end{proof}

    \begin{lemma}[Broadness passes to finer tubes]
        \label{lem-broad-transfer}
        Let $3\leq k\leq n-1$, $0<\rho<1$ and $C_1\geq 2$ with $C_1\rho\leq 1$. Let $\mathcal{S}$ be a finite nonempty family of tubes such that for every $(k-1)$-plane $V$,
        \[\#\{S\in\mathcal{S}: \angle(S,V)\leq C_1\rho\}\leq \frac{1}{2k}\,\#\mathcal{S}.\]
        For each $S\in\mathcal{S}$, let $\T_S$ be a finite family of tubes with $\angle(T,S)\leq \rho$ for every $T\in\T_S$, and suppose that $\mu \leq \#\T_S\leq 2\mu$ for every $S$, for some $\mu\geq 1$. Let $\T := \bigsqcup_{S}\T_S$. Then at least $(5/12)^{k-1}(\#\T)^k$ of the $k$-tuples $(T_1,\ldots,T_k)\in\T^k$ satisfy
        \begin{equation}
            \label{eq-transverse-tuple}
            \angle\big(\dir(T_i), \mathrm{span}\{\dir(T_1),\ldots,\dir(T_{i-1})\}\big) \geq C_1\rho/2 \qquad\text{for } 2\leq i\leq k,
        \end{equation}
        and every such $k$-tuple satisfies $|\dir(T_1)\wedge\cdots\wedge\dir(T_k)| \geq (C_1\rho/4)^{k-1}$. In particular (taking $\T_S = \{S\}$), the same holds for $k$-tuples in $\mathcal{S}^k$.
    \end{lemma}
    \begin{proof}
        Choose $T_1\in\T$ arbitrarily. Suppose $T_1,\ldots,T_{i-1}$ have been chosen, where $2\leq i\leq k$. Their directions span a subspace $W$ of dimension at most $k-1$. Let $V\supseteq W$ be a $(k-1)$-plane. At most $\#\mathcal{S}/(2k)$ of the tubes $S\in\mathcal{S}$ satisfy $\angle(S,V)\leq C_1\rho$. If $\angle(S,V) > C_1\rho$ and $T\in\T_S$, then
        \[\angle(T,W)\geq\angle(T,V)\geq \angle(S,V) - \angle(T,S) \geq C_1\rho - \rho \geq C_1\rho/2.\]
        The number of such $T$ is at least $(1-\frac{1}{2k})\#\mathcal{S}\cdot\mu \geq \frac{1}{2}(1-\frac{1}{2k})\#\T \geq \frac{5}{12}\#\T$, since $k\geq 3$. Multiplying the number of choices at each step gives at least $(5/12)^{k-1}(\#\T)^k$ such $k$-tuples. Finally, for unit vectors, $|v_1\wedge\cdots\wedge v_k| = \prod_{i=2}^k \mathrm{dist}(v_i,\mathrm{span}\{v_1,\ldots,v_{i-1}\})$, and each factor is at least $\sin(C_1\rho/2)\geq C_1\rho/4$.
    \end{proof}

    The following consequence of Zahl's $k$-linear estimate gives a volume bound for families of tubes pointing in $\delta$-separated directions that are $k$-broad at every point.

    \begin{lemma}[A consequence of Zahl's $k$-linear estimate]
        \label{lem-zahl-corollary}
        Let $2\leq k\leq n$, $d := d_{\mathrm{lin}}(k) := (n^2+n+k^2-k)/(2n)$ and $\varepsilon>0$. There is $C_\varepsilon>0$ such that the following holds. Let $0<\delta<1$, $0<\beta\leq 1$ and $0<c\leq 1$, and let $\T$ be a family of $\delta$-tubes pointing in $\delta$-separated directions, with shading $Y$, such that for every $p\in E_\T$,
        \begin{equation}
            \label{eq-zahl-corollary-hyp}
            \#\big\{(T_1,\ldots,T_k)\in\T(p)^k: |\dir(T_1)\wedge\cdots\wedge\dir(T_k)|>\beta\big\}\geq c\,\mu(p)^k.
        \end{equation}
        Then
        \[|E_\T|\geq C_\varepsilon\,c^{d/k}\,\beta\,\delta^{n-d+\varepsilon}\Big(\sum_{T\in\T}|Y(T)|\Big)^{d}\Big(\sum_{T\in\T}|T|\Big)^{-\frac{n(d-1)}{n-1}}.\]
        In particular, if $|Y(T)|\geq\lambda|T|$ for every $T\in\T$, then
        \[|E_\T|\geq C_\varepsilon\,c^{d/k}\,\beta\,\lambda^{d}\,\delta^{n-d+\varepsilon}\Big(\sum_{T\in\T}|T|\Big)^{\frac{n-d}{n-1}}.\]
    \end{lemma}
    Since $\lambda\gtrsim\delta$ (a nonempty shading contains a $\delta$-cube) and $\sum_{T\in\T}|T|\lesssim 1$, the last bound also holds, with a different constant, if $d$ is replaced on the right-hand side by any $d_1\leq d_{\mathrm{lin}}(k)$.
    \begin{proof}
        By Zahl's direction-separated multilinear Kakeya estimate~\cite[Theorem~1.3]{ZAHL2021107596}, for every $\varepsilon'>0$ there is $C_{\varepsilon'}$ such that every family $\T_0$ of $\delta$-tubes pointing in $\delta$-separated directions satisfies
        \begin{equation}
        \label{eq-k-linear-est}
        \begin{split}
        \Big\|\Big(\sum_{T_1,\ldots,T_k\in\T_0}\chi_{T_1}\cdots\chi_{T_k}\,|\dir(T_1)\wedge\cdots\wedge\dir(T_k)|^{k/d}\Big)^{1/k}\Big\|_{L^{d'}}&\\
        \leq C_{\varepsilon'}\,\delta^{1-\frac{n}{d}-\varepsilon'}\Big(\sum_{T\in\T_0}|T|\Big)^{\frac{n(d-1)}{(n-1)d}}&,
        \end{split}
        \end{equation}
        with $d = d_{\mathrm{lin}}(k)$. For $p\in\R^n$ let
        \[F(p) := \Big(\sum_{T_1,\ldots,T_k\in\T}\chi_{Y(T_1)}(p)\cdots\chi_{Y(T_k)}(p)\,|\dir(T_1)\wedge\cdots\wedge\dir(T_k)|^{k/d}\Big)^{1/k}.\]
        For $p\in E_\T$, the $k$-tuples in $(\ref{eq-zahl-corollary-hyp})$ each contribute more than $\beta^{k/d}$ to the sum, so $F(p)\geq c^{1/k}\beta^{1/d}\mu(p)$. Since $\int\mu = \sum_{T}|Y(T)|$ and $F$ is supported in $E_\T$, H\"older's inequality and $(\ref{eq-k-linear-est})$ with $\varepsilon' = \varepsilon/d$ (using $\chi_{Y(T_i)}\leq\chi_{T_i}$) give
        \[c^{1/k}\beta^{1/d}\sum_{T\in\T}|Y(T)|\leq\int F\leq|E_\T|^{1/d}\,\|F\|_{L^{d'}}\leq C_{\varepsilon/d}\,|E_\T|^{1/d}\,\delta^{1-\frac nd-\frac{\varepsilon}{d}}\Big(\sum_{T\in\T}|T|\Big)^{\frac{n(d-1)}{(n-1)d}}.\]
        Raising this to the power $d$ gives the first bound, and the second follows from $\sum_T|Y(T)|\geq\lambda\sum_T|T|$ and $d-\frac{n(d-1)}{n-1} = \frac{n-d}{n-1}$.
    \end{proof}

    Extremal families need not point in $\delta$-separated directions. The next lemma extends Lemma~\ref{lem-zahl-corollary} to families in which at most $\delta^{-t}$ tubes are essentially parallel to any given tube, by passing to a random subfamily that points in $\delta$-separated directions.

    \begin{lemma}[Broad volume estimate]
        \label{lem-broad-volume}
        Let $3\leq k\leq n-1$, $d := d_{\mathrm{lin}}(k) := (n^2+n+k^2-k)/(2n)$ and $\varepsilon'>0$. There is $C = C(n,k,\varepsilon')$ such that the following holds. Let $0<\delta<1$, $t\geq 0$ and $0<b\leq 1$ with $b^{k-1}\geq 4\delta$. Let $\T$ be a family of essentially distinct $\delta$-tubes, at most $\delta^{-t}$ of which are essentially parallel to any given $\delta$-tube, with shading $Y$, and let $\mu(p) := \#\T(p)$. Suppose that for every $p\in E_\T$, at least $(5/12)^{k-1}\mu(p)^k$ of the $k$-tuples $(T_1,\ldots,T_k)\in\T(p)^k$ satisfy $|\dir(T_1)\wedge\cdots\wedge\dir(T_k)|\geq b^{k-1}$. Then
        \[\sum_{T\in\T}|Y(T)|\leq C\,\delta^{-(k+3)t}\,b^{-\frac{(k+3)(k-1)}{d}}\,\delta^{1-\frac{n}{d}-\varepsilon'}\,|E_\T|^{1/d}.\]
    \end{lemma}
    \begin{proof}
        Call two tubes of $\T$ \emph{conflicting} if they are essentially parallel. Each tube conflicts with at most $\delta^{-t}$ other tubes. We will first choose a subfamily of tubes, in which no tubes are conflicting, using probabilistic methods. An important feature of the probabilistic method used is that a transverse $k$-tuple of tubes survives with an explicit probability.
        
        Set $q := \delta^t/2$. Select each tube of $\T$ independently with probability $q$, and let $\T_\omega$ be the set of selected tubes T such that no tube conflicting with T was selected. (For example, if two conflicting tubes are both selected, both are thrown away.) Then, since $\T$ is essentially distinct, tubes in $\T_\omega$ point in $\delta$-separated directions. Further, if $T_1, \ldots, T_k \in \T$ are any $k$ tubes such that $|\dir(T_1)\wedge\cdots\wedge\dir(T_k)|\geq b^{k-1}\geq 4\delta$, then $|\dir(T_i)-\dir(T_j)|>2\delta$ for $i\neq j$, so no tube conflicts with two of the $T_i$, and the $T_i$ do not conflict with each other. Hence the events that $T_i\in\T_\omega$ depend on disjoint sets of tubes, and
        \[\mathbb{P}\big(T_1,\ldots,T_k\in\T_\omega\big) \geq q^k(1-q)^{k\delta^{-t}}\geq  (q/2)^k.\]
        For $p\in E_\T$, let
        \[F_\omega(p) := \Big(\sum_{T_1,\ldots,T_k\in\T_\omega}\chi_{Y(T_1)}(p)\cdots\chi_{Y(T_k)}(p)\,|\dir(T_1)\wedge\cdots\wedge\dir(T_k)|^{k/d}\Big)^{1/k}.\]
        Then $F_\omega(p)\leq\mu(p)$, and by the hypothesis, $\mathbb{E}_\omega\big[F_\omega(p)^k\big]\geq a_0\,\mu(p)^k$, where $a_0 := (q/2)^k(5/12)^{k-1}b^{(k-1)k/d}$. For a random variable $X\in[0,1]$ and $s>0$, $\mathbb{E}[X^k]\leq s^k+\mathbb{P}(X\geq s)$. Applying this to $X = F_\omega(p)/\mu(p)$ with $s := (a_0/2)^{1/k}$ gives 
        \[\mathbb{P}\left(F_\omega(p)\geq s\mu(p)\right)\geq \frac{a_0}{2}.\] Integrating against $\mu$, there is a choice of $\omega$ such that the set $G := \{p\in E_\T: F_\omega(p)\geq s\mu(p)\}$ satisfies \[\int_G\mu\geq\frac{a_0}{2}\int\mu = \frac{a_0}{2}\sum_{T\in\T}|Y(T)|.\] 
        The set $G$ is a union of $\delta$-cubes.
        Give each $T\in\T_\omega$ the shading $Y_\omega(T) := Y(T)\cap G$, and let $\mu_\omega(p) := \#\{T\in\T_\omega: p\in Y_\omega(T)\}$. For $p\in G$, the sum defining $F_\omega(p)^k$ has at most $\mu_\omega(p)^k$ terms, each at most $1$, so $\mu_\omega(p)\geq F_\omega(p)\geq s\mu(p)$. Let $\beta := (s^k/2)^{d/k}$. The $k$-tuples with $|\dir(T_1)\wedge\cdots\wedge\dir(T_k)|\leq\beta$ contribute at most 
        \[\beta^{k/d}\mu_\omega(p)^k = \frac{s^k}{2}\mu_\omega(p)^k\] 
        to $F_\omega(p)^k\geq s^k\mu_\omega(p)^k$, so at least $\frac{s^k}{2}\mu_\omega(p)^k$ of the $k$-tuples in $\T_\omega(p)^k$ have $|\dir(T_1)\wedge\cdots\wedge\dir(T_k)|>\beta$. Moreover \[\sum_{T\in\T_\omega}|Y_\omega(T)| = \int_G\mu_\omega\geq s\int_G\mu\geq\frac{sa_0}{2}\sum_{T\in\T}|Y(T)|,\] and $\sum_{T\in\T_\omega}|T|\lesssim 1.$ Lemma~\ref{lem-zahl-corollary}, applied to $(\T_\omega,Y_\omega)$ with $c = s^k/2$ and this value of $\beta$, gives
        \[|E_\T|\geq|G|\gtrsim C_{\varepsilon'}\,s^{2d}\,\delta^{n-d+\varepsilon'}\Big(\frac{sa_0}{2}\sum_{T\in\T}|Y(T)|\Big)^d.\]
        Hence $\sum_T|Y(T)|\lesssim C_{\varepsilon'}^{-1/d}s^{-3}a_0^{-1}\delta^{1-\frac nd-\varepsilon'}|E_\T|^{1/d}$, and  the conclusion is reached by observing from the definitions of $s$ and $a_0$ that $s^{-3}a_0^{-1}\lesssim a_0^{-(k+3)/k}\lesssim\delta^{-(k+3)t}b^{-(k+3)(k-1)/d}$.
    \end{proof}

    \begin{proposition}
        \label{prop-sigma-bound}
        For every $3\leq k\leq n-1$,
        \[\sigma_n\leq n-\min\big(d_{\mathrm{pl}}(k),\,d_{\mathrm{lin}}(k)\big),\qquad d_{\mathrm{lin}}(k) := \frac{n^2+k^2+n-k}{2n},\]
        where $d_{\mathrm{pl}}(k)$ is as in Theorem~\ref{prop-final-planebrush-est}.
    \end{proposition}
    \begin{proof}
        Fix $k$, write $d_{\mathrm{pl}} := d_{\mathrm{pl}}(k)$, $d_{\mathrm{lin}} := d_{\mathrm{lin}}(k)$, $d := \min(d_{\mathrm{pl}},d_{\mathrm{lin}})$, $f := f(k)$ which is the exponent of the shading density in Theorem~\ref{prop-final-planebrush-est}, and let $0<\varepsilon<(n-2)^{-2}$. 
        Let $\tau := \sqrt{\varepsilon}$, so that $\varepsilon\leq\tau\leq 1-\varepsilon$ and $\tau(k-1)<1$. 
        Let $\eta$ and $\delta_0$ be given by Proposition~\ref{prop-wz-self-similar} for this $\varepsilon$. By decreasing $\eta$, assume $\eta\leq\varepsilon$. Let $\delta\in(0,\delta_0]$ be such that an $\eta$-extremal pair $(\T,Y)_\delta$ exists. 
        By~\cite[Lemma~3.1]{Wang2022StickyKS}, $\delta$ can be taken arbitrarily small. Throughout, $c>0$ denotes a constant that may depend on $n$, $k$ and $\varepsilon$ but not on $\delta$, and $o(1)$ denotes a quantity tending to $0$ as $\delta\to 0$ with $\varepsilon$ fixed. It will be shown that
        \begin{equation}
            \label{eq-sigma-bound}
            \sigma_n\leq n-d+C_n\sqrt{\varepsilon}+o(1),
        \end{equation}
        with $C_n$ depending only on $n$. Letting $\delta\to 0$ and then $\varepsilon\to 0$ proves the proposition.

        Let $\rho := \delta^\tau\in[\delta^{1-\varepsilon},\delta^\varepsilon]$, and let $(\T',Y')_\delta$ and $(\tilde\T,\tilde Y)_\rho$ be given by Proposition~\ref{prop-wz-self-similar}. Since $(\T',Y')_\delta$ is a refinement of an $\eta$-extremal pair,
        \begin{equation}
            \label{eq-fine-mass}
            M_{\mathrm{f}} := \sum_{T\in\T'}|Y'(T)|\geq\delta^{2\eta}\qquad\text{and}\qquad|E_{\T'}|\leq|E_\T|\leq\delta^{\sigma_n-\eta}
        \end{equation}
        for $\delta$ small. Let $a := \delta^{6\varepsilon}$, and apply Lemma~\ref{lem-k-linear-or-k-plany} to $(\tilde\T,\tilde Y)$ at scale $\rho$, with $C_1 = 4$ and threshold $a$.

        \emph{Case (a): mostly plany.} By Proposition~\ref{prop-wz-self-similar}(ii), $\sum_{\tilde T}|Z'(\tilde T)|\geq\frac a2\sum_{\tilde T}|\tilde Y(\tilde T)|\geq\frac a2\rho^{\varepsilon}$. Join two tubes of $\tilde\T$ by an edge if they are essentially parallel. Each tube has at most $\rho^{-\varepsilon}$ neighbours, so $\tilde\T$ splits into at most $2\rho^{-\varepsilon}$ classes, each consisting of tubes pointing in $\rho$-separated directions. One class carries at least $\frac a4\rho^{2\varepsilon}$ of the mass of $Z'$. By dyadic pigeonholing, it has a subfamily $\mathcal{S}$ and a number $\lambda$ with $\lambda|S|\leq|Z'(S)|<2\lambda|S|$ for $S\in\mathcal{S}$ and $\sum_{S\in\mathcal{S}}|Z'(S)|\geq c\,a\rho^{2\varepsilon}/\log(1/\delta)$. Since $\sum_{S\in\mathcal{S}}|S|\lesssim 1$, this gives $\lambda\geq c\,a\rho^{2\varepsilon}/\log(1/\delta)\geq\delta^{9\varepsilon}$ for $\delta$ small. The family $(\mathcal{S},Z')$ is $(k-1)$-plany with constant $4$ at scale $\rho$, its tubes point in $\rho$-separated directions, and $Z'(S)\subseteq S$. Applying Theorem~\ref{prop-final-planebrush-est} at scale $\rho$, with $\varepsilon$ in place of its $\varepsilon$, gives
        \begin{align*}|E_{\tilde\T}|\geq\Big|\bigcup_{S\in\mathcal{S}}Z'(S)\Big|&\geq c\,\lambda^{f}\rho^{n-d_{\mathrm{pl}}+\varepsilon}\sum_{S\in\mathcal{S}}|S|\geq c\,\lambda^{f}\rho^{n-d_{\mathrm{pl}}+\varepsilon}\sum_{S\in\mathcal{S}}|Z'(S)|\\ &\geq c\,\delta^{(9f+8)\varepsilon}\rho^{n-d_{\mathrm{pl}}+\varepsilon}/\log(1/\delta).\end{align*}
        Combined with $|E_{\tilde\T}|\leq\rho^{\sigma_n-\varepsilon}$ from Proposition~\ref{prop-wz-self-similar}(ii), and recalling $\rho = \delta^\tau$,
        \[\tau(\sigma_n-\varepsilon)\leq(9f+8)\varepsilon+\tau(n-d_{\mathrm{pl}}+\varepsilon)+o(1),\]
        that is, $\sigma_n\leq n-d_{\mathrm{pl}}+2\varepsilon+(9f+8)\sqrt\varepsilon+o(1)$. Since $f\leq n$, this implies $(\ref{eq-sigma-bound})$.

        \emph{Case (b): mostly broad.} Let $B$ be the union of the broad $\rho$-cubes. For $p\in\R^n$, let $\nu(p) := \#\{\tilde T\in\tilde\T: p\in E_{\T'[\tilde T]}\}$ and $\mu(p) := \#\T'(p)$. By Proposition~\ref{prop-wz-self-similar}(i) and (iv), $\mu = \mu_{\mathrm{fine}}\,\nu$. Every $Y'(T)$ with $T\in\T'[\tilde T]$ lies in $E_{\T'[\tilde T]}\subseteq\tilde Y(\tilde T)$, so by (iv) and (v),
        \[M_{\mathrm{f}} = \sum_{\tilde T}\sum_{Q\subseteq\tilde Y(\tilde T)}\mu_{\mathrm{fine}}\,|Q\cap E_{\T'[\tilde T]}| = \mu_{\mathrm{fine}}\,w\,\#\{(\tilde T,Q): Q\subseteq\tilde Y(\tilde T)\},\]
        and the same computation restricted to $B$ shows that $\int_B\mu$ is the fraction of pairs $(\tilde T, Q)$ with $Q$ broad, times $M_{\mathrm{f}}$. In case (b) this fraction is at least $\frac12$, so $\int_B\mu\geq\frac12M_{\mathrm{f}}$. Let $\nu_0 := M_{\mathrm{f}}/(4\mu_{\mathrm{fine}}|E_{\T'}|)$ and $G := B\cap\{\nu\geq\nu_0\}$. Since $\int_{\{\nu<\nu_0\}}\mu\leq\mu_{\mathrm{fine}}\nu_0|E_{\T'}| = \frac14M_{\mathrm{f}}$,
        \begin{equation}
            \label{eq-good-mass}
            \int_G\mu\geq\tfrac14M_{\mathrm{f}}\geq\tfrac14\delta^{2\eta}.
        \end{equation}
        By $(\ref{eq-fine-mass})$ and the bound on $\mu_{\mathrm{fine}}$ in (iv), $\nu_0\geq\frac14\delta^{3\eta}(\delta/\rho)^{\sigma_n+\varepsilon}\delta^{-\sigma_n}\geq\frac14\delta^{3\eta+\varepsilon}\rho^{-\sigma_n}$. If $Q$ is a $\rho$-cube, then $\#\{\tilde T: Q\subseteq\tilde Y(\tilde T)\}\leq\rho^{-\sigma_n-\varepsilon}$ by (iii), so for $p\in G\cap Q$,
        \[\nu(p)\geq\nu_0\geq\tfrac14\delta^{3\eta+2\varepsilon}\,\#\{\tilde T: Q\subseteq\tilde Y(\tilde T)\}\geq\tfrac14\delta^{5\varepsilon}\,\#\{\tilde T: Q\subseteq\tilde Y(\tilde T)\}.\]

        Let $p\in G$, let $Q$ be the $\rho$-cube containing $p$, and let $V$ be a $(k-1)$-plane. If $T\in\T'(p)$ and $T\in\T'[\tilde T]$, then $p\in E_{\T'[\tilde T]}\subseteq\tilde Y(\tilde T)$, so $Q\subseteq\tilde Y(\tilde T)$, and $\angle(T,\tilde T)\leq\rho/2$. Hence if $\angle(T,V)\leq\rho$, then $\angle(\tilde T,V)\leq 4\rho$. Each $\tilde T$ contributes at most $\mu_{\mathrm{fine}}$ tubes to $\T'(p)$, and $Q$ is broad, so
        \[\#\{T\in\T'(p): \angle(T,V)\leq\rho\}\leq\mu_{\mathrm{fine}}\,a\,\#\{\tilde T: Q\subseteq\tilde Y(\tilde T)\}\leq 4\delta^{\varepsilon}\mu_{\mathrm{fine}}\,\nu(p) = 4\delta^\varepsilon\mu(p)\leq\frac{1}{2k}\mu(p)\]
        for $\delta$ small. By Lemma~\ref{lem-broad-transfer}, with $\mathcal{S} = \T'(p)$, $\T_S = \{S\}$, $\rho/2$ in place of $\rho$ and $C_1 = 2$, at least $(5/12)^{k-1}\mu(p)^k$ of the $k$-tuples in $\T'(p)^k$ have $|\dir(T_1)\wedge\cdots\wedge\dir(T_k)|\geq b^{k-1}$, where $b := \rho/\pi$. Since $\tau(k-1)<1$, $b^{k-1}\geq 4\delta$ for $\delta$ small.

        Let $Y''(T) := Y'(T)\cap G$, a union of $\delta$-cubes, so that $\#\{T: p\in Y''(T)\} = \mu(p)$ for $p\in G$. By property (b) of $(\T,Y)_\delta$ at scale $\delta$, at most $C\delta^{-\eta}\leq\delta^{-2\eta}$ tubes of $\T'$ are essentially parallel to any given $\delta$-tube. Lemma~\ref{lem-broad-volume} with $t = 2\eta$ and $\varepsilon' = \varepsilon$, together with $(\ref{eq-good-mass})$, gives
        \[\tfrac14\delta^{2\eta}\leq\sum_{T}|Y''(T)|\leq c^{-1}\delta^{-2(k+3)\eta}\rho^{-\frac{(k+3)(k-1)}{d_{\mathrm{lin}}}}\delta^{1-\frac{n}{d_{\mathrm{lin}}}-\varepsilon}|E_{\T}|^{1/d_{\mathrm{lin}}}.\]
        Hence $|E_\T|\geq c\,\delta^{n-d_{\mathrm{lin}}+C_n\varepsilon}\rho^{(k+3)(k-1)}$, and combined with $|E_\T|\leq\delta^{\sigma_n-\eta}$ this gives $\sigma_n\leq n-d_{\mathrm{lin}}+(k+3)(k-1)\sqrt\varepsilon+C_n\varepsilon+o(1)$, which implies $(\ref{eq-sigma-bound})$.
    \end{proof}

    \begin{proposition}
        \label{prop-sticky-estimate}
        Let
        \[d = \max_{3\leq k\leq n-1}\min \left( d_{\mathrm{pl}}(k), \frac{n^2+k^2+n-k}{2n}\right),\]
        with $d_{\mathrm{pl}}(k)$ as in Theorem~\ref{prop-final-planebrush-est}. For every $\kappa > 0$ there are $\varepsilon_0>0$ and $\delta_0>0$ such that, whenever $0<\varepsilon\leq\varepsilon_0$, $0<\delta\leq\delta_0$, and $\T$ is a sticky family of $\delta$-tubes with shading $Y$ satisfying $|Y(T)| \geq \delta^\varepsilon |T|$ for all $T \in \T$,
         \[|E_\T| \geq \delta^{n-d+\kappa}. \]
    \end{proposition}
    \begin{proof}
        By Proposition~\ref{prop-sigma-bound}, $\sigma_n\leq n-d$. Let $\kappa>0$. Since $N(s,t)\to\sigma_n$ as $(s,t)\to(0,0)$, there are $s_0,t_0\in(0,1)$ with $N(s_0,t_0)\leq n-d+\kappa/2$. Let $\varepsilon_0 := s_0/3$. By the definition of $N$ and its monotonicity, there is $\delta_0>0$ such that $M(3\varepsilon,t_0,\delta)\geq M(s_0,t_0,\delta)\geq\delta^{n-d+\kappa}$ for all $\varepsilon\leq\varepsilon_0$ and $\delta\leq\delta_0$.

        Now let $\T$ be a sticky family of $\delta$-tubes with $|Y(T)|\geq\delta^\varepsilon|T|$. After a rigid motion and a dilation by a constant factor, and after passing to the tubes whose directions lie in a fixed cap (which contains at least a constant fraction of them), it may be assumed that the tubes of $\T$ are $\delta$-tubes in the sense above, up to constant factors in $\delta$ that only affect the constants. Then $\T$ satisfies (a), (b) and (c) in Definition $\ref{def:epsilon-extremal}$ with $s = 3\varepsilon$ and $t = t_0$, for $\delta$ small:
        \begin{itemize}
            \item (a) holds because the directions are $\delta$-separated.
            \item For (b), given $\delta\leq\rho\leq 1$, let $\theta := \rho/C$ for a suitable constant $C$. The $\theta$-parents of $\T$ point in $\theta$-separated directions, and each $T\in\T$ lies in its $\theta$-parent. So the $\rho$-tubes with the same coaxial lines as the $\theta$-parents cover $\T$, and $O(1)\leq\delta^{-t_0}$ of them are essentially parallel to a common $\rho$-tube. For $\rho<C\delta$, the tubes of $\T$ themselves can be used.
            \item (c) holds because $\sum_{T}|Y(T)|\geq\delta^{\varepsilon}\,\#\T\,|T|\gtrsim\delta^{2\varepsilon}\geq\delta^{3\varepsilon}$.
        \end{itemize}
        Hence $|E_\T|\geq M(3\varepsilon,t_0,\delta)\geq\delta^{n-d+\kappa}$.
    \end{proof}

    The proof of Theorem \ref{cor-sticky-estimate} follows immediately from this result.
    \begin{proof}[Proof of Theorem \ref{cor-sticky-estimate}]
        By Proposition~\ref{prop-sticky-estimate}, it suffices to find an integer $3 \leq k \leq n-1$ with
        \[\min\left( d_{\mathrm{pl}}(k), \frac{n^2+k^2+n-k}{2n}\right) \geq 0.60376707287\, n + \mathcal{O}(1).\]
        Recall from Theorem~\ref{prop-final-planebrush-est} that, for $k \geq 12$,
        \[d_{\mathrm{pl}}(k) = \frac{2n+2-(\sqrt{2}-1)k}{3} - \frac{n-2-2(\sqrt{2}-1)k}{3(3(n-k)+2)}.\]
        First consider the simpler expression
        \[\min \left( \frac{2n+2-(\sqrt{2}-1)k}{3}, \frac{n^2+k^2+n-k}{2n}\right).\]
        Its maximum over real $k$ is attained where the two expressions agree, at
        \[k = \frac{1-\sqrt{2} + \sqrt{6-2\sqrt{2}}}{3}n + \mathcal{O}(1) \approx 0.45555915723 \,n + \mathcal{O}(1),\]
        and there both expressions equal $0.60376707287\, n + \mathcal{O}(1)$. Let $k_0$ be the nearest integer to this value of $k$; rounding changes the minimum by only $\mathcal{O}(1)$. For large $n$ we have $k_0 \geq 12$ and $n - k_0 \geq n/2$, so the explicit formula applies and the correction term in $d_{\mathrm{pl}}(k_0)$ is at most
        \[\frac{n}{9(n-k_0)} \leq \frac{2}{9}.\]
        Hence
        \[\min\left( d_{\mathrm{pl}}(k_0), \frac{n^2+k_0^2+n-k_0}{2n}\right) \geq 0.60376707287\, n + \mathcal{O}(1),\]
        and Proposition~\ref{prop-sticky-estimate} holds with
        \[d = 0.60376707287 \, n +\mathcal{O}(1).\]
    \end{proof}
    
    \subsection{Relation to the sticky Besicovitch sets of Wang and Zahl}
    \label{subsec:wang-zahl}
    Wang and Zahl~\cite{Wang2022StickyKS} defined \emph{sticky Besicovitch sets} using packing dimension. This subsection explains how sticky Besicovitch sets in the sense of Wang and Zahl give rise to sticky families of $\delta$-tubes. A version of this result is proved in~\cite{Wang2022StickyKS}.
    \begin{definition}
        A \emph{sticky Besicovitch set} is a  compact set $K \subseteq \R^n$ obeying the following hypotheses: there is a set of lines $L$ in $\R^n$ such that 
        \begin{enumerate}[i.]
            \item $L$ has a line in every direction, 
            \item $K\cap \ell$  contains a unit line segment for every $\ell \in L$, and
            \item  $L$ has packing dimension $n-1$.
        \end{enumerate}
    \end{definition}
    Lines are parametrized, and $\dist$ and $m_{\mathcal{L}}$ are defined, as in Subsection~\ref{subsec:sticky-kakeya}.
    
    \begin{proposition}
        For any $\varepsilon > 0$, 
        if $K$ is a sticky Besicovitch set in $\R^n$,
        then there exists $\delta_0 > 0$ such that 
        \begin{enumerate}[(1)]
            \item for any $\delta < \delta_0$, there is a family of $\delta$-tubes $\T$ which is sticky; and
            \item if $\dim_H(K) = d$, then there is a shading $Y(T)$ satisfying $|Y(T)|\gtrsim  \big(\log(1/\delta)\big)^{-2} |T|$ for every $T\in \T$, such that
            \[|E_\T| \lesssim \delta^{n-d-\varepsilon}.\]

        \end{enumerate}
    \end{proposition}
    \begin{proof}
    Let $\varepsilon > 0$ be given and 
    let $K$ be a sticky Besicovitch set. 
    Since the packing dimension of $L$ coincides with its upper modified box dimension, $L$ can be written as a countable union of closed sets $\{L_i\}$ such that 
    \[\overline{\dim_M}(L_i) \leq n-1+\frac{\varepsilon}{4} \]
    for all $i$. There exists an index $i$ for which $\nu\big(\dir(L_i)\big)>0$. From the upper Minkowski dimension of $L_i$, there exists $\delta_1 > 0$ such that for all $\delta < \delta_1$,
    \[m_{\mathcal{L}}(\mathcal{N}_\delta(L_i)) \leq \delta^{(2n-2)-(n-1)-\varepsilon/2}.\]
    A direction $e\in\mathbb{S}^{n-1}$ is said to be \emph{over-represented at scale} $\delta$ if
    \[|\{\textbf{a}\in \R^{n-1}: (\textbf{a},e)\in \mathcal{N}_\delta(L_i)\}| > \delta^{n-1-\varepsilon}. \]
    By Fubini's theorem, it follows that
    \[\nu\big(\{e\in \mathbb{S}^{n-1}: e \text{ is over-represented at scale } \delta \}\big) \leq \delta^{\varepsilon/2}.\]
    By making $\delta_1$ smaller if necessary, it can be ensured that
    \[\sum_{\text{ dyadic }\delta: \delta < \delta_1} \nu\big(\{e\in \mathbb{S}^{n-1}: e \text{ is over-represented at scale } \delta \}\big) \leq \frac{1}{2} \nu(\dir(L_i)).\]
    Now define a subset $L' \subseteq L_i$ by
    \[L' = L_i \setminus \bigcup_{\text{ dyadic }\delta: \delta < \delta_0} \dir^{-1}\big(\{e\in \mathbb{S}^{n-1}: e \text{ is over-represented at scale } \delta \}\big). \]
    By the definition of Hausdorff dimension, there exists a number $\delta_0 > 0$ (which may be taken smaller than $\delta_1$ above) and a cover $\mathcal{Q}$ of $K$, where \[\mathcal{Q} = \bigcup_{\text{dyadic} \delta: \delta < \delta_0} \mathcal{Q}_\delta,\]
    where $\mathcal{Q}_\delta$ is a collection of $\delta$-cubes and
    \[\#\mathcal{Q}_\delta \lesssim \delta^{-d-\varepsilon/2}.\]
    For each $\ell \in L'$, since $|K\cap \ell| \geq 1$,
    there exists at least one $\delta$ such that
    \begin{equation}
    \label{eq-wang-zahl-1}
    \bigg|\ell \cap \big(\bigcup_{Q\in \mathcal{Q}_\delta} Q\big) \bigg| \geq \frac{1}{100\big(\log(1/\delta)\big)^2}.
    \end{equation}
    For each dyadic number $\delta$, let $L_\delta \subseteq L'$ be the set of lines for which $(\ref{eq-wang-zahl-1})$ holds. Then there exists a dyadic scale $\delta < \delta_0$ for which $\nu(L_\delta)> \nu(L')/100 \big(\log(1/\delta)\big)^2$. Fix this value of $\delta$ and rename $L_\delta$ as $L''$.
    Now let $\mathcal{E}$ be a $\delta$-separated subset of $\dir(L'')$. Then
    \[\#\mathcal{E} \geq \frac{\nu(L')\delta^{1-n}}{100\big(\log(1/\delta)\big)^2}.\]
    Define $\T$ to be a family of $\delta$-tubes pointing in directions in $\mathcal{E}$, with centers and shading $Y(T)$ chosen so that, using $(\ref{eq-wang-zahl-1})$,
    \[E_\T = \bigcup_{Q \in \mathcal{Q}_\delta}Q\]
    and $|Y(T)| \gtrsim |T|/\big(\log(1/\delta)\big)^2$ for every $T\in \T$. 
    Since this implies that
    \[|E_\T| \lesssim \delta^{n-d-\varepsilon}, \]
    conclusion (2) in the statement of the proposition is shown.
    It now remains to show that $\T$ is sticky.
    Let $\rho$ be a number such that $\delta < \rho < 1$. 
    
    If $\rho \geq \delta_0$: since $\delta_0$ may be considered an admissible constant, if $\T_\rho$ is any set of maximal $\rho$-separated tubes, it is trivial that $\#\T_\rho \approx \rho^{1-n}$ and $\#\T[T_\rho] \approx (\delta/\rho)^{1-n}.$

    If $\rho < \delta_0$: cover $L''$ greedily with $\rho$-balls, so that $\rho/3$-balls at the same centers would be disjoint.
    Define $\T_\rho$ to be the set of tubes pointing in the directions of the lines in the cover and with centers inside the set on the left-hand-side of $(\ref{eq-wang-zahl-1})$. From this definition, it is clear that every $T\in \T$ is contained in some $T_\rho \in \T_\rho$.

    Since $\delta_0 < \delta_1$, and since $\dir(L'')$ has no directions over-represented at any scale less than $\delta_1$, there are at most $\lesssim \delta^{-\varepsilon}$ tubes $T_\rho \in \T_\rho$ that point in the same direction. 
    To make the directions $\rho$-separated, choose any one tube out of the possible $\delta^{-\varepsilon}$ tubes in the cover pointing in a given direction, and remove the other tubes pointing in that direction from $\T_\rho$. Then $\#\T_\rho \approx \rho^{1-n}$. 
    To ensure that $\T$ is still covered by $\T_\rho$, take a refinement $\T' \subset \T$, losing at most a factor of $\rho^{-\varepsilon} \leq \delta^{-\varepsilon}$ in the number of tubes. Then $\T'$ will be a sticky set, and it can be easily seen that conclusion (2) still holds for this refinement.
    \end{proof}
    \begin{corollary}
    \label{cor-sticky-besicovitch}
        The Hausdorff dimension of every sticky Besicovitch set in $\R^n$ is at least $0.60376707287 \,n + \mathcal{O}(1)$. 
    \end{corollary}
    \begin{proof}
        Let $d$ be as in Theorem~\ref{cor-sticky-estimate}, and suppose that $K$ is a sticky Besicovitch set with $\dim_H(K) = d' < d$. Let $\kappa := (d-d')/3$, let $\varepsilon_0$ and $\delta_0$ be given by Theorem~\ref{cor-sticky-estimate}, and let $\varepsilon := \min(\varepsilon_0, \kappa)$. The previous proposition, with this $\varepsilon$, gives arbitrarily small $\delta$ and sticky families $\T$ of $\delta$-tubes with shadings satisfying $|Y(T)| \gtrsim (\log(1/\delta))^{-2}|T| \geq \delta^{\varepsilon}|T|$ and $|E_\T| \lesssim \delta^{n-d'-\varepsilon} \leq \delta^{n-d'-\kappa}$. For $\delta$ small, this contradicts the lower bound $\delta^{n-d+\kappa}$ from Theorem~\ref{cor-sticky-estimate}, since $n-d'-\kappa > n-d+\kappa$. Alternatively, the corollary follows from~\cite[Proposition~2.1]{Wang2022StickyKS} and Proposition~\ref{prop-sigma-bound}.
    \end{proof}

    \section{X-ray estimates and Minkowski dimension}
    \label{sec:xray}
        Katz, {\L}aba and Tao~\cite{Katz-Laba-Tao} showed that an x-ray estimate, that is, a Kakeya maximal estimate with an improved dependence on the number of tubes pointing in each direction, reduces Minkowski dimension bounds for Besicovitch sets to the sticky case. This section records an x-ray estimate at every dimension below
        \[d_x := \max_{2\leq k \leq n}\min\Big\{\frac{n^2 + n+ k^2 -k}{2n},\, n-k+2\Big\} = (2-\sqrt{2})\,n + \mathcal{O}(1),\]
        obtained by interpolating the Kakeya maximal estimate of Hickman, Rogers and Zhang~\cite{Hickman2019ImprovedBF} and Zahl~\cite{ZAHL2021107596} with the x-ray estimate of {\L}aba and Tao~\cite{Laba1999xray} (Proposition~\ref{prop-x-ray}). It then shows that an x-ray estimate at a dimension $d$ below the exponent of Theorem~\ref{cor-sticky-estimate} gives a Minkowski dimension bound strictly larger than $d$ (Corollary~\ref{cor-Minkowski}).

    \begin{definition}[X-ray estimate]
        Let $d\geq 1$ and $0\leq\alpha,\gamma\leq 1$, and write $d' := d/(d-1)$, with $d' = \infty$ if $d = 1$. An \emph{x-ray bound with exponents $(\alpha,\gamma)$ holds at dimension $d$} if for every $\varepsilon>0$ there is $C_\varepsilon>0$ such that the following holds for every $0<\delta<1$ and $m\geq 1$. If $\mathcal{E}\subseteq\mathbb{S}^{n-1}$ is a set of $\delta$-separated directions and $\T$ is a family of $\delta$-tubes such that
        \begin{enumerate}[i.]
            \item no tube in $\T$ is contained in the $10$-fold dilate of another one, and
            \item every tube in $\T$ points in a direction in $\mathcal{E}$, and at most $m$ tubes in $\T$ point in any given direction,
        \end{enumerate}
        then
        \[\Big\|\sum_{T\in\T}\chi_T\Big\|_{L^{d'}}\leq C_\varepsilon\,\delta^{1-\frac{n}{d}-\varepsilon}\,m^{\alpha}\,\big(\delta^{n-1}\#\mathcal{E}\big)^{\gamma}.\]
        An \emph{x-ray estimate holds at dimension $d$} if an x-ray bound holds at dimension $d$ with some exponents $\alpha,\gamma<1$.
    \end{definition}
    Splitting $\T$ into $m$ families pointing in $\delta$-separated directions shows that a Kakeya maximal estimate at dimension $d$ gives an x-ray bound at dimension $d$ with $\alpha = 1$; the point of an x-ray estimate is that $\alpha<1$.

    \begin{lemma}[Interpolation]
        \label{lem-x-ray-interpolation}
        Suppose that x-ray bounds hold at dimensions $d_a, d_b\geq 1$ with exponents $(\alpha_a,\gamma_a)$ and $(\alpha_b,\gamma_b)$. Let $0\leq\vartheta\leq 1$ and define $d_0$ by $\frac{1}{d_0} = \frac{1-\vartheta}{d_a}+\frac{\vartheta}{d_b}$. Then an x-ray bound holds at dimension $d_0$ with exponents $\big((1-\vartheta)\alpha_a+\vartheta\alpha_b,\,(1-\vartheta)\gamma_a+\vartheta\gamma_b\big)$.
    \end{lemma}
    \begin{proof}
        Since $\frac{1}{d_0'} = 1-\frac{1}{d_0} = \frac{1-\vartheta}{d_a'}+\frac{\vartheta}{d_b'}$, H\"older's inequality gives $\|f\|_{L^{d_0'}}\leq\|f\|_{L^{d_a'}}^{1-\vartheta}\|f\|_{L^{d_b'}}^{\vartheta}$ for $f = \sum_{T\in\T}\chi_T$. Inserting the two bounds, the powers of $m$ and of $\delta^{n-1}\#\mathcal{E}$ are as claimed, and the power of $\delta$ is $(1-\vartheta)(1-\frac{n}{d_a})+\vartheta(1-\frac{n}{d_b}) = 1-\frac{n}{d_0}$, up to $\delta^{-\varepsilon}$.
    \end{proof}

    \begin{lemma}[Three x-ray bounds]
        \label{lem-x-ray-inputs}
        \begin{enumerate}[(a)]
            \item An x-ray bound holds at dimension $1$ with exponents $(1,1)$.
            \item An x-ray bound holds at dimension $d_x$ with exponents $\big(1,\frac{d_x-1}{d_x}\big)$.
            \item An x-ray bound holds at dimension $p := \frac{n+2}{2}$ with exponents
            \[\alpha_{\mathrm{LT}} := 1-\frac{1}{2(n+2)},\qquad\gamma_{\mathrm{LT}} := 1-\frac{n}{(n-1)(n+2)}.\]
        \end{enumerate}
    \end{lemma}
    \begin{proof}
        (a) At most $m\#\mathcal{E}$ tubes pass through any point, so $\|\sum_T\chi_T\|_{L^\infty}\leq m\#\mathcal{E} = \delta^{1-n}m(\delta^{n-1}\#\mathcal{E})$.

        (b) By~\cite{Hickman2019ImprovedBF} and~\cite[Theorem~1.5]{ZAHL2021107596}, the Kakeya maximal estimate (Definition~\ref{def-maximal-estimate}) holds in $\R^n$ at dimension $d_x$. Split $\T$ into $m$ subfamilies $\T_1,\ldots,\T_m$, where $\T_j$ contains the $j$-th tube in each direction (in some ordering). Each $\T_j$ points in $\delta$-separated directions and has $\#\T_j\leq\#\mathcal{E}$, so the maximal estimate and the triangle inequality give $\|\sum_T\chi_T\|_{L^{d_x'}}\leq m\,C_\varepsilon\delta^{1-\frac{n}{d_x}-\varepsilon}(\delta^{n-1}\#\mathcal{E})^{(d_x-1)/d_x}$.

        (c) Let $q := \frac{(n-1)(n+2)}{n}$ and $r := 2(n+2)$. By~\cite[Theorem~1.2]{Laba1999xray}, if the tubes of $\T$ point in directions in a $\delta$-separated set, at most $m$ in each direction, and tubes with the same direction have $\delta$-separated coaxial lines, then
        \[\Big\|\sum_{T\in\T}\chi_T\Big\|_{L^{p'}}\lesssim\delta^{1-\frac{n}{p}-\varepsilon}\,m^{\frac1q-\frac1r}\,\big(\delta^{n-1}\#\T\big)^{1/q'}.\]
        Condition i implies the separation of coaxial lines after splitting $\T$ into $O(1)$ subfamilies (and rotating each of $O(1)$ caps of directions into the range considered in~\cite{Laba1999xray}). Since $\#\T\leq m\#\mathcal{E}$, the right-hand side is at most a constant times $\delta^{1-\frac np-\varepsilon}m^{1-\frac1r}(\delta^{n-1}\#\mathcal{E})^{1/q'}$, and $1-\frac1r = \alpha_{\mathrm{LT}}$, $\frac{1}{q'} = \gamma_{\mathrm{LT}}$.
    \end{proof}

    \begin{proposition}
        \label{prop-x-ray}
        An x-ray estimate holds in $\R^n$ at every dimension $1<d_0<d_x$.
    \end{proposition}
    \begin{proof}
        Let $p = \frac{n+2}{2}$. If $d_0\leq p$, apply Lemma~\ref{lem-x-ray-interpolation} to the bounds (c) and (a) of Lemma~\ref{lem-x-ray-inputs}, with $\frac{1}{d_0} = \frac{1-\vartheta}{p}+\vartheta$. Then $0\leq\vartheta<1$ because $d_0>1$, so the resulting exponents $(1-\vartheta)\alpha_{\mathrm{LT}}+\vartheta$ and $(1-\vartheta)\gamma_{\mathrm{LT}}+\vartheta$ are less than $1$. If $p<d_0<d_x$, apply Lemma~\ref{lem-x-ray-interpolation} to the bounds (b) and (c), with $\frac{1}{d_0} = \frac{1-\vartheta}{d_x}+\frac{\vartheta}{p}$. Then $0<\vartheta<1$, so the exponents $(1-\vartheta)+\vartheta\alpha_{\mathrm{LT}} = 1-\frac{\vartheta}{2(n+2)}$ and $(1-\vartheta)\frac{d_x-1}{d_x}+\vartheta\gamma_{\mathrm{LT}}$ are less than $1$.
    \end{proof}

    The next proposition is a version of~\cite[Proposition~3.3]{Katz-Laba-Tao}. It produces a family of tubes satisfying conditions (a), (b) and (c) of the definition of $M(s,t,\delta)$ in Subsection~\ref{subsec:sticky-kakeya}; tubes, lines and the measure $m_{\mathcal{L}}$ are as defined there.
    \begin{proposition}
        \label{prop-x-ray-implies-sticky}
        Suppose that an x-ray estimate holds at dimension $1<d<n$ with exponents $\alpha,\gamma<1$, and let $C := \frac{4(d+2)}{d(1-\alpha)}$. Let $\varepsilon>0$, and let $B$ be a Besicovitch set with $\overline{\dim}_M(B)<d+\varepsilon$. Then for every sufficiently small $\delta>0$ there is a family $(\T,Y)_\delta$ satisfying conditions (a), (b) and (c) of the definition of $M(s,t,\delta)$ with $s = \varepsilon$ and $t = C\varepsilon$, and
        \[|E_\T|\leq\delta^{n-d-\varepsilon}.\]
    \end{proposition}
    \begin{proof}
        After a translation and a dilation by a factor depending only on $n$, assume that $B\subseteq[-\frac{1}{2n},\frac{1}{2n}]^n$ and that $B$ contains a segment of length $c_n>0$ in every direction. Let $\Omega := \{e\in\mathbb{S}^{n-1}: \angle(e,e_n)\leq\pi/3\}$, and for each $e\in\Omega$ fix a line $\ell_e = (\textbf{a}_e, e)$ such that $\ell_e\cap B$ contains a segment of length $c_n$. Let $L := \{\ell_e: e\in\Omega\}$. Since $\overline{\dim}_M(B)<d+\varepsilon$, there is $r_B>0$ such that $|\mathcal{N}_r(B)|\leq r^{n-d-\varepsilon}$ for $0<r\leq r_B$. All implicit constants below depend only on $n$, $d$, $\alpha$, $\gamma$ and the constants in the x-ray estimate, but not on $\delta$.

        Let $C_1 := \frac{d+2}{d(1-\alpha)} = C/4$. For $A\geq 1$ and $0<r<1$, call a direction $e\in\Omega$ \emph{sticky at scale $r$} if
        \begin{equation}
        \label{eq-sticky-direction-condition}
        \big|\{\textbf{a}\in\R^{n-1}: (\textbf{a},e)\in\mathcal{N}_r(L)\}\big|\leq A\,r^{n-1-C_1\varepsilon},
        \end{equation}
        where $\mathcal{N}_r(L)$ is the $r$-neighborhood of $L$ with respect to the distance between lines, and let $D_r\subseteq\Omega$ be the set of directions that are not sticky at scale $r$.

        \emph{Step 1: for every $0<r\leq r_B/C'$, $\nu(D_r)\lesssim r^{c_1\varepsilon}$, where $c_1 := \frac{1}{d(1-\gamma)}$.} Here $C'$ is a large constant depending only on $n$. Let $\rho := C'r$ and let $\mathcal{E}_1$ be a maximal $\rho$-separated subset of $D_r$, so that $\nu(D_r)\lesssim\rho^{n-1}\#\mathcal{E}_1$. For $e\in\mathcal{E}_1$, the set in $(\ref{eq-sticky-direction-condition})$ has measure greater than $Ar^{n-1-C_1\varepsilon}$, so it contains $m\gtrsim Ar^{-C_1\varepsilon}$ points $\textbf{a}_1,\ldots,\textbf{a}_m$ with pairwise distances at least $20\rho$; fix exactly $m := \lceil cAr^{-C_1\varepsilon}\rceil$ of them. Each line $(\textbf{a}_i,e)$ is within distance $r$ of a line of $L$, so it contains a segment of length $c_n$ lying in $\mathcal{N}_{C''r}(B)$. Let $\T_0$ be the family of $\rho$-neighborhoods of these segments, over all $e\in\mathcal{E}_1$ and $1\leq i\leq m$; if $C'$ is large enough, $\T_0$ satisfies conditions i and ii of the definition of an x-ray bound at scale $\rho$, with $\mathcal{E} = \mathcal{E}_1$, and $E_{\T_0}\subseteq\mathcal{N}_{2\rho}(B)$. Write $X := \rho^{n-1}\#\mathcal{E}_1$. Since $\sum_{T\in\T_0}|T|\sim mX$, H\"older's inequality and the x-ray estimate give
        \[mX\lesssim|E_{\T_0}|^{1/d}\Big\|\sum_{T\in\T_0}\chi_T\Big\|_{L^{d'}}\lesssim|E_{\T_0}|^{1/d}\,\rho^{1-\frac nd-\varepsilon}\,m^{\alpha}X^{\gamma},\]
        so that, using $|E_{\T_0}|\leq|\mathcal{N}_{2\rho}(B)|\lesssim\rho^{n-d-\varepsilon}$,
        \[X^{d(1-\gamma)}\lesssim\rho^{-(d+1)\varepsilon}\,m^{-d(1-\alpha)}\lesssim A^{-d(1-\alpha)}\,r^{-(d+1)\varepsilon+dC_1(1-\alpha)\varepsilon} = A^{-d(1-\alpha)}\,r^{\varepsilon}.\]
        Hence $\nu(D_r)\lesssim X\lesssim r^{c_1\varepsilon}$, since $A\geq 1$.

        \emph{Step 2: choice of the sticky directions.} By Step 1, $\sum\nu(D_r)\lesssim\sum_{r\leq r_1}r^{c_1\varepsilon}\lesssim r_1^{c_1\varepsilon}/(c_1\varepsilon)$, where the sums are over dyadic $r$. Choose $r_1\leq r_B/C'$, depending on $\varepsilon$ but not on $\delta$, so small that this is at most $\nu(\Omega)/2$. Then choose $A := C'''r_1^{1-n}$ with $C'''$ large, so that every direction is sticky at every scale $r\geq r_1/8$, because the set in $(\ref{eq-sticky-direction-condition})$ always has measure $\lesssim 1$. Let
        \[\mathcal{E} := \Omega\setminus\bigcup_{\substack{r \text{ dyadic}\\ \delta\leq r\leq r_1}}D_r,\qquad\text{so that }\nu(\mathcal{E})\geq\tfrac12\nu(\Omega).\]
        Every $e\in\mathcal{E}$ is sticky at every dyadic scale $r\geq\delta$.

        \emph{Step 3: the family of tubes.} Let $\mathcal{E}_\delta$ be a maximal $\delta$-separated subset of $\mathcal{E}$, and for $e\in\mathcal{E}_\delta$ let $T_e$ be the $\delta$-tube with coaxial line $\ell_e$, with shading $Y(T_e)$ the union of the $\delta$-cubes that meet the segment of length $c_n$ in $\ell_e\cap B$. Let $\T := \{T_e: e\in\mathcal{E}_\delta\}$. Then $E_\T\subseteq\mathcal{N}_{\sqrt n\,\delta}(B)$, so $|E_\T|\leq\delta^{n-d-\varepsilon}$ for $\delta$ small, after decreasing $\varepsilon$ slightly at the start. Condition (a) holds because the directions are $\delta$-separated, and (c) holds because $\sum_T|Y(T)|\gtrsim\delta^{n-1}\#\mathcal{E}_\delta\gtrsim\nu(\mathcal{E})\gtrsim 1\geq\delta^{\varepsilon}$ for $\delta$ small.

        For (b), let $\delta\leq\rho\leq 1$. As in~\cite[proof of Proposition~2.1]{Wang2022StickyKS}, choose greedily a set $\tilde L\subseteq\{\ell_e: e\in\mathcal{E}_\delta\}$ such that the $\rho/2$-neighborhoods of the lines in $\tilde L$ cover $\{\ell_e: e\in\mathcal{E}_\delta\}$ and their $\rho/6$-neighborhoods are disjoint. The $\rho$-tubes with coaxial lines in $\tilde L$ cover $\T$. Suppose that $N$ of them are essentially parallel to a common $\rho$-tube with direction $v$. Their coaxial lines $\tilde\ell_1,\ldots,\tilde\ell_N$ have directions within $\rho$ of $v$, so the disjoint sets $\mathcal{N}_{\rho/6}(\tilde\ell_i)$ lie in $\mathcal{N}_{\rho/6}(L)\cap\dir^{-1}(\mathcal{N}_{2\rho}(\{v\}))$, and
        \[m_{\mathcal{L}}\Big(\mathcal{N}_{\rho/6}(L)\cap\dir^{-1}(\mathcal{N}_{2\rho}(\{v\}))\Big)\gtrsim N\rho^{2(n-1)}.\]
        By Fubini's theorem, there is $v'$ with $|v'-v|\leq 2\rho$ and $|\{\textbf{a}: (\textbf{a},v')\in\mathcal{N}_{\rho/6}(L)\}|\gtrsim N\rho^{n-1}$. Let $e := \dir(\tilde\ell_1)\in\mathcal{E}$, so that $|e-v'|\leq 3\rho$. If $(\textbf{a},v')\in\mathcal{N}_{\rho/6}(L)$, then $(\textbf{a},e)\in\mathcal{N}_{4\rho}(L)$, so $|\{\textbf{a}: (\textbf{a},e)\in\mathcal{N}_{r}(L)\}|\gtrsim N\rho^{n-1}$ for the dyadic $r\in[4\rho,8\rho)$. Since $e$ is sticky at scale $r$, $N\lesssim A\rho^{-C_1\varepsilon}\leq\delta^{-2C_1\varepsilon}$ for $\delta$ small. Hence (b) holds with $t = 2C_1\varepsilon\leq C\varepsilon$.
    \end{proof}

    \begin{corollary}
    \label{cor-Minkowski}
        Let $d_{\mathrm{st}} := \max_{3\leq k\leq n-1}\min\big(d_{\mathrm{pl}}(k), \frac{n^2+k^2+n-k}{2n}\big) = 0.60376707287\,n + \mathcal{O}(1)$ be the exponent of Proposition~\ref{prop-sticky-estimate}. If an x-ray estimate holds at a dimension $1<d<d_{\mathrm{st}}$, then there is $\varepsilon_0>0$, depending on $n$, $d$ and the x-ray estimate, such that every Besicovitch set in $\R^n$ has upper Minkowski dimension at least $d+\varepsilon_0$.
    \end{corollary}
    \begin{proof}
        Let $\kappa := (d_{\mathrm{st}}-d)/4$, and let $C$ be as in Proposition~\ref{prop-x-ray-implies-sticky}. By Proposition~\ref{prop-sigma-bound}, $\sigma_n\leq n-d_{\mathrm{st}}$, so there are $s_0,t_0\in(0,1)$ with $N(s_0,t_0)\leq n-d_{\mathrm{st}}+\kappa$, and hence $M(s_0,t_0,\delta)\geq\delta^{n-d_{\mathrm{st}}+2\kappa}$ for all sufficiently small $\delta$. Let $\varepsilon_0 := \min(s_0, t_0/C, \kappa)$, and suppose that $B$ is a Besicovitch set with $\overline{\dim}_M(B)<d+\varepsilon_0$. Proposition~\ref{prop-x-ray-implies-sticky} with $\varepsilon = \varepsilon_0$ gives, for arbitrarily small $\delta$, a family satisfying (a), (b), (c) with $s = \varepsilon_0\leq s_0$ and $t = C\varepsilon_0\leq t_0$, so that, by monotonicity of $M$,
        \[\delta^{n-d_{\mathrm{st}}+2\kappa}\leq M(s_0,t_0,\delta)\leq M(\varepsilon_0,C\varepsilon_0,\delta)\leq|E_\T|\leq\delta^{n-d-\varepsilon_0}.\]
        Hence $d_{\mathrm{st}}-d\leq 2\kappa+\varepsilon_0\leq 3\kappa<d_{\mathrm{st}}-d$, a contradiction.
    \end{proof}
    Proposition~\ref{prop-x-ray} gives x-ray estimates only at dimensions $d_0<d_x$, and Corollary~\ref{cor-Minkowski} then gives the lower bound $d_0+\varepsilon_0(d_0)$, where $\varepsilon_0(d_0)$ may tend to $0$ as $d_0\to d_x$. This does not improve on the Minkowski dimension bound $d_x$, which already follows from the Kakeya maximal estimate of~\cite{Hickman2019ImprovedBF,ZAHL2021107596}. Since $d_x = (2-\sqrt2)n+\mathcal{O}(1)<d_{\mathrm{st}}$ for large $n$, an x-ray estimate at dimension $d_x$ itself would, by Corollary~\ref{cor-Minkowski}, show that every Besicovitch set in $\R^n$ has upper Minkowski dimension at least $d_x+\varepsilon_0$.

\bmhead{Acknowledgments}
I am grateful to my PhD advisor, Professor Nets Katz, for suggesting this problem and for supervising the project. I would like to thank Professor Joshua Zahl for many helpful discussions, and for pointing out that the x-ray estimate in Section~\ref{sec:xray} can be derived as a corollary of his $k$-linear Kakeya theorem.


\begin{thebibliography}{14}
\providecommand{\natexlab}[1]{#1}
\providecommand{\url}[1]{\texttt{#1}}
\expandafter\ifx\csname urlstyle\endcsname\relax
  \providecommand{\doi}[1]{doi: #1}\else
  \providecommand{\doi}{doi: \begingroup \urlstyle{rm}\Url}\fi

\bibitem[Bourgain(1991)]{Bourgain1991}
J.~Bourgain.
\newblock {B}esicovitch type maximal operators and applications to {F}ourier
  analysis.
\newblock \emph{Geometric and Functional Analysis}, 1\penalty0 (2):\penalty0
  147--187, 1991.
\newblock URL \url{http://eudml.org/doc/58112}.

\bibitem[C\'ordoba(1977)]{Cordoba1977}
Antonio C\'ordoba.
\newblock The {K}akeya maximal function and the spherical summation
  multipliers.
\newblock \emph{American Journal of Mathematics}, 99\penalty0 (1):\penalty0
  1--22, 1977.
\newblock URL \url{http://www.jstor.org/stable/2374006}.

\bibitem[Hickman et~al.(2022)Hickman, Rogers, and Zhang]{Hickman2019ImprovedBF}
Jonathan Hickman, Keith~M. Rogers, and Ruixiang Zhang.
\newblock Improved bounds for the {K}akeya maximal conjecture in higher
  dimensions.
\newblock \emph{American Journal of Mathematics}, 144\penalty0 (6):\penalty0
  1511--1560, 2022.

\bibitem[Katz and Tao(2002)]{Katz-Tao}
Nets~Hawk Katz and Terence Tao.
\newblock New bounds for {K}akeya problems.
\newblock \emph{Journal d'Analyse Mathématique}, 87:\penalty0 231--263, 2002.
\newblock \doi{10.1007/BF02868476}.

\bibitem[Katz and Zahl(2021)]{Katz-Zahl-planebrush}
Nets~Hawk Katz and Joshua Zahl.
\newblock A {K}akeya maximal function estimate in four dimensions using
  planebrushes.
\newblock \emph{Revista Matemática Iberoamericana}, 37\penalty0 (1):\penalty0
  317--359, 2021.
\newblock \doi{10.4171/rmi/1219}.

\bibitem[Katz et~al.(2000)Katz, {\L}aba, and Tao]{Katz-Laba-Tao}
Nets~Hawk Katz, Izabella {\L}aba, and Terence Tao.
\newblock An improved bound on the {M}inkowski dimension of {B}esicovitch sets
  in $\mathbb{R}^3$.
\newblock \emph{Annals of Mathematics}, 152\penalty0 (2):\penalty0 383--446,
  2000.
\newblock URL \url{http://www.jstor.org/stable/2661389}.

\bibitem[{\L}aba and Tao(2001)]{Laba1999xray}
Izabella {\L}aba and Terence Tao.
\newblock An {X}-ray transform estimate in $\mathbb{R}^n$.
\newblock \emph{Revista Matemática Iberoamericana}, 17\penalty0 (2):\penalty0
  375--407, 2001.
\newblock \doi{10.4171/RMI/298}.

\bibitem[{OpenAI}(2026)]{OpenAI2026LocalSmoothing}
{OpenAI}.
\newblock Local smoothing for the {E}uclidean wave equation, 2026.
\newblock URL
  \url{https://github.com/openai/math/blob/main/preprints/Local-smoothing-for-the-Euclidean-wave-equation-October-25-2026/paper.pdf}.
\newblock Manuscript dated 25 October 2026.

\bibitem[Tao(1999)]{Tao1999Lectures}
Terence Tao.
\newblock Lecture notes for {M}ath 254{B}.
\newblock University of California, Los Angeles, 1999.
\newblock URL \url{https://www.math.ucla.edu/~tao/254b.1.99s/}.

\bibitem[Wang and Zahl(2022)]{Wang2022StickyKS}
Hong Wang and Joshua Zahl.
\newblock Sticky {K}akeya sets and the sticky {K}akeya conjecture, 2022.
\newblock arXiv:2210.09581.

\bibitem[Wang and Zahl(2025)]{WangZahl2025Kakeya3D}
Hong Wang and Joshua Zahl.
\newblock Volume estimates for unions of convex sets, and the {K}akeya set
  conjecture in three dimensions, 2025.
\newblock arXiv:2502.17655.

\bibitem[Wolff(1995)]{Wolff1995}
Thomas Wolff.
\newblock An improved bound for {K}akeya type maximal functions.
\newblock \emph{Revista Matemática Iberoamericana}, 11\penalty0 (3):\penalty0
  651--674, 1995.
\newblock URL \url{http://eudml.org/doc/39495}.

\bibitem[Wolff(1998)]{Wolff-x-ray}
Thomas Wolff.
\newblock A mixed norm estimate for the {X}-ray transform.
\newblock \emph{Revista Matemática Iberoamericana}, 14\penalty0 (3):\penalty0
  561--600, 1998.
\newblock URL \url{http://eudml.org/doc/39560}.

\bibitem[Zahl(2021)]{ZAHL2021107596}
Joshua Zahl.
\newblock New {K}akeya estimates using {G}romov's algebraic lemma.
\newblock \emph{Advances in Mathematics}, 380:\penalty0 107596, 2021.
\newblock \doi{10.1016/j.aim.2021.107596}.

\end{thebibliography}
\end{document}